\documentclass[12pt]{amsart}
\usepackage{amsmath,amsfonts,amsthm,amssymb,amscd}
\usepackage[mathscr]{eucal}
\allowdisplaybreaks
\usepackage{graphicx}
\usepackage{times}%%%
\usepackage[curve,matrix,arrow]{xy}
\xyoption{poly}
\advance\textheight\topskip
\newcommand{\algW}{\boldsymbol{\mathcal{W}}}%%%{\mathfrak{W}}
\newcommand{\qgr}{\mathfrak{G}}

\newcommand{\vac}{|0\rangle}
\newcommand{\muvac}[1]{|{\mu_{#1}}\rangle}

\newcommand{\logid}{\mathbb{L}}

\newcommand{\bref}[1]{\textbf{\ref{#1}}}

\newcommand{\cZ}{\mathfrak{Z}}

\newcommand{\im}{\mathop{\mathrm{im}}\nolimits}
\renewcommand{\ker}{\mathop{\mathrm{ker}}\nolimits}
\newcommand{\injection}{\rightarrowtail}

\newcommand{\surjection}{\twoheadrightarrow}

\newcommand{\embedding}{\injection}%%%{\hookrightarrow}

\renewcommand{\geq}{\,{\geqslant}\,}
\renewcommand{\leq}{\,{\leqslant}\,}

\renewcommand{\le}{\,{\leqslant}\,}

\newcommand{\abp}{z_1\,{:}\,z_2}
\newcommand{\HomSL}{\mathrm{Hom}_{\UresSL2}}
\newcommand{\Hom}{\mathrm{Hom}_{\overline{\mathscr{U}}_{\q}}}
\newcommand{\HomC}{\mathrm{Hom}_{\mathbb{C}}}
\newcommand{\justHom}{\mathrm{Hom}}
\newcommand{\HomW}[1]{\mathrm{Hom}_{\algW(#1)}}

\newcommand{\Ext}{\mathrm{Ext}_{\rule{0pt}{6.5pt}%
    \overline{\mathscr{U}}_{\q}}^1}
\newcommand{\EXT}{\mathrm{Ext}^\bullet}
\newcommand{\Extn}{\mathrm{Ext}_{\rule{0pt}{9.5pt}%
    \overline{\mathscr{U}}_{\q}}^n}

\newcommand{\Yext}[1]{\mathrm{Ext}_{\rule{0pt}{9.5pt}%
    \overline{\mathscr{U}}_{\q}}^#1}

\newcommand{\justExt}{\mathrm{Ext}}
\newcommand{\morX}{\mathsf{X}}
\newcommand{\ob}{\mathrm{Ob}}
\newcommand{\bP}{\mathbf{P}}

\newcommand{\Rep}{\mathrm{Rep}} %% a category of representations
\newcommand{\qK}{\boldsymbol{\mathrm{K}}} %% The Kronecker quiver
\newcommand{\Kron}{\boldsymbol{\mathrm{K}}} %% The Kronecker quiver
\newcommand{\Aone}{\mathrm{A^{(1)}_1}}

\newcommand{\q}{\mathfrak{q}}

\newcommand{\tensor}{\otimes}

\newcommand{\lc}{\rep{C}_p}

\newcommand{\Mmod}{\mathbb{M}}

\newcommand{\fusion}{%
  \mathop{{\otimes}\kern-7pt\raisebox{.6pt}{%
      \mbox{\footnotesize${\bullet}$}}}}

\newcommand{\ket}[1]{|#1\rangle}

\newcommand{\UresSL}[1]{\overline{\mathscr{U}}_{\q} s\ell(#1)}

\newcommand{\mfrac}[2]{\mbox{\small$\displaystyle\frac{#1}{#2}$}}
\newcommand{\ffrac}[2]{\mbox{\footnotesize$\displaystyle\frac{#1}{#2}$}}
\newcommand{\half}{%
  \mathchoice{\ffrac{1}{2}\,}{\frac{1}{2}}{\frac{1}{2}}{\frac{1}{2}}}
\newcommand{\fhalf}{\ffrac{1}{2}}

\newcommand{\bs}[1]{\boldsymbol{#1}}
\newcommand{\modPr}{\mathscr{P}}
\newcommand{\WrepP}{\mathfrak{P}}

\newcommand{\modM}{\mathscr{M}}
\newcommand{\modN}{\mathscr{N}} %<-------------- aNy module

\newcommand{\modO}{\mathscr{O}}
\newcommand{\repX}{\mathscr{X}}
\newcommand{\WrepX}{\mathfrak{X}}

\newcommand{\repReg}{\mathscr{R}}%<---------------- left Regular module
\newcommand{\modW}{\mathscr{W}}

\newcommand{\modI}{\mathscr{I}}

\newcommand{\Verma}{\mathscr{V}}
\newcommand{\CVerma}{\bar{\mathscr{V}}}

\newcommand{\repPi}{\repX^{-}}

\newcommand{\rep}{\mathscr}  %% representations
\newcommand{\Wrep}{\mathfrak}  %% W-alg representations
\newcommand{\setI}{\mathrm{Ind}} %% a set of indecomposable modules

\newcommand{\FunF}{\mathcal{F}} %% a functor
\newcommand{\FunG}{\mathcal{G}} %% a functor
\newcommand{\dd}{\partial}
\newcommand{\SLiiZ}{SL(2,\oZ)}
\newcommand{\oC}{\mathbb{C}}
\newcommand{\oP}{\mathbb{P}}
\newcommand{\oN}{\mathbb{N}}

\newcommand{\oZ}{\mathbb{Z}}

\newcommand{\one}{\boldsymbol{1}}%{1\kern-4pt 1}

\newcommand{\toppr}{\mathsf{b}}
\newcommand{\botpr}{\mathsf{a}}
\newcommand{\leftpr}{\mathsf{x}}
\newcommand{\rightpr}{\mathsf{y}}

\newcommand{\ribbon}{{\boldsymbol{v}}}

\newcommand{\cas}{\boldsymbol{C}}

\numberwithin{equation}{section}
\makeatletter
\@addtoreset{equation}{section}
\@addtoreset{subsubsection}{section}

\def\@secnumfont{\bfseries}
\def\subsubsection{\@startsection{subsubsection}{3}%
  \z@{.5\linespacing\@plus.7\linespacing}{-.5em}%
  {\normalfont\bfseries}}
\def\paragraph{\@startsection{paragraph}{4}%
  \z@\z@{-\fontdimen2\font}%
  \normalfont\bfseries}
\def\subparagraph{\@startsection{subparagraph}{5}%
  \z@\z@{-\fontdimen2\font}%
  \normalfont\bfseries}

\makeatother

\swapnumbers
\newtheorem{Thm}[subsection]{Theorem}
\newtheorem{Conj}[subsection]{Conjecture}
\newtheorem{Thmm}[subsection]{Classification theorem for $\lc$}

\newtheorem{Lemma}[subsection]{Lemma}
\newtheorem{lemma}[subsubsection]{Lemma}

\newtheorem{prop}[subsubsection]{Proposition}

\theoremstyle{definition}
\newtheorem{Dfn}[subsection]{Definition}

\newtheorem{rem}[subsubsection]{Remark}

\begin{document}
\title[Representation category in logarithmic CFT models in
quantum-group terms]{%
  \vspace*{-3\baselineskip}
  \mbox{}\hfill\texttt{\small\lowercase{math}.QA/0512621}
  \\[2\baselineskip]
  Kazhdan--Lusztig correspondence for the representation category of
  the triplet $W$-algebra in logarithmic CFT}

\author[Feigin]{B.L.~Feigin}%

\address{\mbox{}\kern-\parindent  
  blf: Landau Institute for Theoretical Physics
  \hfill\mbox{}\linebreak 
  \texttt{feigin@mccme.ru}}

\author[Gainutdinov]{A.M.~Gainutdinov}%

\address{\mbox{}\kern-\parindent  
  amg: Physics Department, Moscow State University
  \hfill\mbox{}\linebreak 
  \texttt{azot@mccme.ru}}

\address{\mbox{}\kern-\parindent  
  ams: Lebedev Physics Institute
  \hfill\mbox{}\linebreak 
  \texttt{asemikha@td.lpi.ru}}

\author[Semikhatov]{A.M.~Semikhatov}

\address{\mbox{}\kern-\parindent  
  iyt: Lebedev Physics Institute
  \hfill\mbox{}\linebreak 
  \texttt{tipunin@td.lpi.ru}}

\author[Tipunin]{I.Yu.~Tipunin}

\begin{abstract}
  To study the representation category of the triplet $W$-algebra
  $\algW(p)$ that is the symmetry of the $(1,p)$ logarithmic conformal
  field theory model, we propose the equivalent category $\lc$ of
  finite-dimensional representations of the restricted quantum group
  $\UresSL2$ at $\q=e^{\frac{i\pi}{p}}$.  We fully describe the
  category $\lc$ by classifying all indecomposable representations.
  These are exhausted by projective modules and three series of
  representations that are essentially described by indecomposable
  representations of the Kronecker quiver.  The equivalence of the
  $\algW(p)$- and $\UresSL2$-representation categories is conjectured
  for all $p\geq2$ and proved for $p=2$, the implications including
  the identifications of the quantum-group center with the logarithmic
  conformal field theory center and of the universal $R$-matrix with
  the braiding matrix.
\end{abstract}

\maketitle

\thispagestyle{empty}

\setcounter{tocdepth}{2}%3

%% \flushcolumns

%% \def\contentsname{{}}

%% \vspace*{-36pt}

%% \noindent\mbox{}\kern-16pt\parbox{1.05\textwidth}{
%%   \begin{multicols}{2}
%%     \renewcommand{\textbf}{\relax}
%%     {\footnotesize
%%       \tableofcontents}
%%   \end{multicols}
%% }

\vspace*{-24pt}

\begin{scriptsize}
  \renewcommand{\textbf}{\relax}
  \addtolength{\baselineskip}{-12pt}
  \tableofcontents
\end{scriptsize}

\section{Introduction}

\noindent
\textbf{1.0.} \ The Kazhdan--Lusztig correspondence~\cite{[KL]}
between a vertex-operator algebra (a conformal field theory) and a
quantum group, although functorial, is not necessarily an equivalence
(nor, strictly speaking, does it necessarily involve the
representation category of a \textit{quantum group}; this often has to
be replaced by a quasitensor category~\cite{[Fink],[T]} that is not
the representation category of any quantum group).  ``Deviations''
from the equivalence typically occur in \textit{rational} conformal
field theories~\cite{[MS],[FbZ],[BK],[fuRs4],[fuRs8]}.  A class of
conformal field theory models where the Kazhdan--Lusztig
correspondence is an equivalence is given by the $(1,p)$ logarithmic
models.  This remarkable property may be attributed to the fact that
the $(1,p)$ models are not logarithmic \textit{extensions} of minimal
conformal field theory models, simply because the latter are
nonexistent.  In the $(1,p)$ model, the representation category is
nonsemisimple, and establishing the Kazhdan--Lusztig correspondence
does not therefore require a special ``semisimplification'' of the
representation category of the quantum group.  For recent advances in
the study of vertex-operator algebras with nonsemisimple
representation categories, see~\cite{[My3],[HLZ]} and the references
therein and in~\cite{[jF]}.

In the $(1,p)$ logarithmic models, the vertex-operator algebra is the
triplet $W$-algebra $\algW(p)$ ``based on the $\phi_{3,1}$
operator''~\cite{[K-first]}, and the Kazhdan--Lusztig-dual quantum
group is $\UresSL2$, the restricted quantum $s\ell(2)$ at
$\q=e^{\frac{i\pi}{p}}$~\cite{[FGST]}.  The $\algW(p)$~algebra and the
$\UresSL2$ quantum group have the same number ($2p$) of irreducible
representations and identical fusion rings (the Grothendieck ring for
the quantum group).  The results in~\cite{[FGST],[FHST]} suggest that
the quasitensor categories of $\algW(p)$- and
$\UresSL2$-representations are in fact equivalent.

In this paper, we prove the equivalence in the most explicit case of
$p=2$.  The proof extends to the general $p$, but requires more
technical details, which will be given elsewhere.  The general
strategy is to construct a bimodule (a right $\UresSL2$-module and a
left $\algW(p)$-module) that decomposes into projective modules from
each side and therefore establishes the equivalence of the two
representation categories.  For $p=2$, this construction is quite
simple because the relevant vertex operators are ``almost local'' with
respect to each other.

By the equivalence, which we assume to hold for any $p$, the study of
a category of $\algW(p)$-representations can be replaced by the study
of the category of finite-dimensional $\UresSL2$-representations.
These are not difficult to classify using the standard
means~\cite{[ARS]}; we show that apart from the projective modules,
all indecomposable representations are enumerated in terms of
indecomposable representations of the Kronecker quiver and can easily
be constructed explicitly; there are four series labeled by
$n\in\oZ_{\geq2}$ and two series labeled by
$(n\in\oZ_{\geq1},z\in\oC\oP^1)$.  In accordance with the equivalence,
this classification also applies to $\algW(p)$-representations;
finding their classification directly would be a much more difficult
problem.

The equivalence of categories leads to other identifications between
``conformal field theory data'' and ``quantum-group data,'' the most
prominent one being the identification of the respective
\textit{centers}.  The center $\cZ_{\mathrm{cft}}$ of a conformal
field theory model\,---\,endomorphisms of the identity functor in the
representation category of the corresponding vertex-operator
algebra~$\mathscr{A}$\,---\,contains ``highly nonlocal'' operators
(the Ishibashi states).  A basis in the center is naturally
constructed in terms of conformal blocks on a torus.  In the
semisimple case, $\cZ_{\mathrm{cft}}$ coincides with the space of
conformal field theory characters, but in logarithmic models, it is
not exhausted by the characters (cf.~\cite{[F]}).  This has been
studied recently within different approaches~\cite{[KeLu],[My3],[FG]}.
The general approach to finding the missing ``extended characters''
via the construction of pseudo-traces for a class of vertex-operator
algebras with nonsemisimple representation category was proposed
in~\cite{[My3]}.  The extended characters, along with the usual ones,
are solutions of certain differential equations~\cite{[FG]}.

When the Kazhdan--Lusztig correspondence establishes the equivalence
between representation categories of a vertex-operator algebra and a
quantum group, $\cZ_{\mathrm{cft}}$ is canonically identified with the
\textit{center~$\cZ$ of the quantum group}.  Moreover, in the $(1,p)$
models at least, this identification between $\cZ_{\mathrm{cft}}$
and~$\cZ$ extends to modular-group representations: the
$\SLiiZ$-representation on $\cZ_{\mathrm{cft}}$ induced by the
$\tau\mapsto-\frac{1}{\tau}$ and $\tau\mapsto\tau+1$ transformations
is equivalent~\cite{[FGST]} to the $\SLiiZ$-action on the quantum
group center~$\cZ$, defined in quite general terms of ribbon braided
tensor categories~\cite{[Lyu],[LM]}.  The actual identification of
$\cZ_{\mathrm{cft}}$ with $\cZ$ yields a remarkable relation
\begin{equation*}
  \ribbon=e^{2i\pi L_0},
\end{equation*}
where~$\ribbon$ is the ribbon element (belonging to the quantum group
center) and $L_0$ is the zero-mode of the energy-momentum tensor.
Furthermore, the $R$-matrix associated with the Kazhdan--Lusztig-dual
quantum group~$\qgr$ can be identified with the braiding matrix in
conformal field theory, its action on vertex operators following from
the construction of~$\qgr$ in terms of screening and contour-removal
operators.\footnote{We note that the Kazhdan--Lusztig-dual quantum
  group $\qgr$ itself is \textit{not} quasitriangular (does not have a
  universal $R$-matrix), at least in the known examples of logarithmic
  $(1,p)$~\cite{[FGST]} and $(p',p)$ models.  But the universal
  $R$-matrix exists for a somewhat larger quantum group, generated by
  the generators of $\qgr$ \textit{and} the square root $k$ of the
  Cartan generator $K$ in $\qgr$, i.e., $k^2=K$.  In conformal field
  theory, the counterpart of the ambiguity in defining $k=\sqrt{K}$ is
  the ambiguity in choosing the direction of contours along which
  operators are permuted.  This leads to \textit{two} braiding
  matrices $B^+$ and $B^-$, well known in (rational) conformal field
  theories~\cite{[MS],[FbZ],[BK],[fuRs4],[fuRs8],[KeLu]}.}

\begin{Dfn}\label{sec:Ures}
  Let $\UresSL2$ be the restricted quantum $s\ell(2)$ at
  $\q=e^{\frac{i\pi}{p}}$, generated by $E$, $F$, and $K$ with the
  relations
  \begin{equation*}
    E^{p}=F^{p}=0,\quad K^{2p}=\one
  \end{equation*}
  and the Hopf-algebra structure given by
  \begin{gather*}
    KEK^{-1}=\q^2E,\quad
    KFK^{-1}=\q^{-2}F,\quad
    [E,F]=\ffrac{K-K^{-1}}{\q-\q^{-1}},\\
    \Delta(E)=\one\otimes E+E\otimes K,\quad
    \Delta(F)=K^{-1}\otimes F+F\otimes\one,\quad
    \Delta(K)=K\otimes K,\\
    \epsilon(E)=\epsilon(F)=0,\quad\epsilon(K)=1,\\
    S(E)=-EK^{-1},\quad  S(F)=-KF,\quad S(K)=K^{-1}.
  \end{gather*}
\end{Dfn}
This quantum group has $2p$ inequivalent irreducible representations,
labeled as $\repX^+_{s}$ and $\repX^-_{s}$, $1\leq s\leq p$.
The module $\repX^{\pm}_{s}$ (with $\dim\repX^{\pm}_{s}=s$) has the
highest weight~$\pm\q^{s-1}$.

Let $\lc$ be the category of finite-dimensional $\oZ$-graded
$\UresSL2$-modules at $\q=e^{\frac{i\pi}{p}}$.  In the standard way,
it is endowed with the structure of a ribbon braided tensor category
(the ribbon element required for this can be found in accordance with
the standard recipe from the universal $R$-matrix found
in~\cite{[FGST]} for a quantum group that contains $\UresSL2$ as a
subalgebra; it turns out that this results in the ribbon element not
for the larger quantum group but just for $\UresSL2$).

\medskip
\addtocounter{subsection}{1}
\newcounter{wsec}
\setcounter{wsec}{\arabic{subsection}}

\noindent
\textbf{\thesubsection. \ The $\algW(p)$ algebra.} For a scalar field
$\varphi(z)$ with the OPE
\begin{equation*}
  \varphi(z)\,\varphi(w)= \log(z-w),
\end{equation*}
the Virasoro algebra with central charge
$c\,{=}\,13\,{-}\,6(p\,{+}\,\frac1p)$ is generated by the
energy-momentum tensor (with normal ordering understood)
\begin{equation*}
  T=\fhalf\,\dd\varphi\,\dd\varphi+\ffrac{\alpha_0}{2}\,\dd^2\varphi,
\end{equation*}
where $\alpha_+=\sqrt{2p}$\,, $\alpha_-=-\sqrt{\frac{2}{p}}$\,, and
$\alpha_0\,{=}\,\alpha_+ \,{+}\, \alpha_-$.  There are the screening
operators
\begin{equation*}
  S_+=\ffrac{1}{2i\pi}\oint\!du\,e^{\alpha_+\varphi(u)},
  \qquad
  S_-=\ffrac{1}{2i\pi}\oint\!du\,e^{\alpha_-\varphi(u)},
\end{equation*}
and the $W$-algebra $\algW(p)$~\cite{[K-first]} is generated by the
three currents
\begin{equation}\label{the-W}
  W^-(z)\!=\!\phi_{3,1}(z)\equiv e^{-\alpha_+\varphi}(z),\quad
  W^0(z)\!=\![S_+,W^-(z)],\quad
  W^+(z)\!=\![S_+,W^0(z)]
\end{equation}
(which are Virasoro primaries of conformal dimension $2p{-}1$).

There are $2p$ irreducible modules of $\algW(p)$, labeled as
$\WrepX^+_{s}$ (``singlet,'' associated with the $\phi_{1,s}$
operator) and $\WrepX^-_{s}$ (``doublet,'' associated with the
$\phi_{2,s}$ operator), $1\leq s\leq p$, with the characters given
by~\cite{[FHST]}
\begin{equation*}%%%\label{eq:characters}
  \begin{aligned}
    \chi_{\WrepX^{+}_{s}}^{\phantom{y}}(q)&=\mfrac{1}{\eta(q)}
    \Bigl(\ffrac{s}{p}\,\theta_{p{-}s,p}(q)
    + 2\,\theta'_{p{-}s,p}(q)\Bigr),\\[4pt]
    \chi_{\WrepX^{-}_{s}}^{\phantom{y}}(q)&= \mfrac{1}{\eta(q)}
    \Bigl(\ffrac{s}{p}\,\theta_{s,p}(q) - 2\,\theta'_{s,p}(q)\Bigr),
  \end{aligned}\qquad
  1\leq s\leq p,
\end{equation*}
where 
\begin{equation*}
  \eta(q)=q^{\frac{1}{24}} \prod_{n=1}^{\infty} (1-q^n),
  \qquad
  \theta_{s,p}(q,z)=\sum_{j\in\oZ
    + \frac{s}{2p}} q^{p j^2} z^j,
\end{equation*}
and we set $\theta_{s,p}(q)=\theta_{s,p}(q,1)$ and
$\theta'_{s,p}(q)=z\frac{\dd}{\dd
  z}\theta_{s,p}(q,z)\!\!\bigm|_{z=1}$.  The fusion for $\algW(2)$ was
directly calculated in~\cite{[GK1]} (also see~\cite{[FHST]} for
$\algW(p)$).

A suitable choice of the representation category for $\algW(p)$ is the
category~$\mathcal{O}$.
\begin{Dfn}
  Let $\Wrep{C}_p$ denote the category $\mathcal{O}$ of
  representations of the $\algW(p)$ algebra, i.e., the category of
  modules in which adjoint eigenspaces of $L_0$ are finite-dimensional
  and the spectrum of $L_0$ is bounded from one side.
\end{Dfn}
This category is a braided ribbon quasitensor category; its
quasitensor structure follows from operator product expansions in
conformal field theory, the braiding is derived from the monodromy of
correlation functions, and the ribbon structure is given by the
operator~$e^{2i\pi L_0}$.

\begin{Thm}\label{Thm:equv}
  The categories $\rep{C}_2$ and $\Wrep{C}_2$ are equivalent as ribbon
  braided quasitensor categories.
\end{Thm}

As noted above, the proof of this fact is greatly simplified by the
existence, for $p=2$, of the symplectic-fermion fields that are
``almost local'' with respect to each other.  In the general case, the
proof requires working out some technical details, but there is little
doubt that the following statement holds.
\begin{Conj}
  For each integer $p\geq2$, the categories $\rep{C}_p$ and
  $\mathfrak{C}_p$ are equivalent as ribbon braided quasitensor
  categories.
\end{Conj}

Assuming the equivalence for all $p$, we can therefore describe the
category of $\algW(p)$-representations by describing the category of
$\UresSL2$-representations, which is a sufficiently simple problem; to
a large degree, it amounts to describing indecomposable
representations of the Kronecker quiver.  We now give a list of
indecomposable modules and then state the result in~\bref{thm:main}.

\medskip
\addtocounter{subsection}{1}
\newcounter{indsec}
\setcounter{indsec}{\arabic{subsection}}

\noindent
\textbf{\thesubsection. \ Indecomposable modules.}
The basic fact is that those spaces $\Ext$ of extensions between the
irreducible representations that are nonzero are two-dimensional
(see~\bref{lem:ext-irrep}); in terms of the bases $\{x^+_i\}$ and
$\{x^-_i\}$, $i=1,2$, chosen in the respective spaces
$\oC^2=\Ext(\repX^{+}_{s},$\linebreak[0]$\repX^{-}_{p-s})$ and
$\oC^2=\Ext(\repX^{-}_{p-s},\repX^{+}_{s})$, we can construct four
families of indecomposable modules as follows.  In the diagrams below,
$\repX_1\stackrel{x}{\longrightarrow}\repX_2$ denotes the extension
by~$x\in\Ext(\repX_1,\repX_2)$, with $\repX_1$ being an irreducible
subquotient and $\repX_2$ an irreducible submodule.
\begin{description}\addtolength{\itemsep}{6pt}
  
\item[$\smash{\boldsymbol{\modW^{{\pm}}_{s}(n)}}$] For $1\leq s \leq
  p{-}1$, $a\,{=}\,\pm$, and integer $n\geq2$, the module
  $\modW^{a}_{s}(n)$ is defined~as

  \includegraphics[bb=1.6in 9.5in 8in 10.6in, clip]{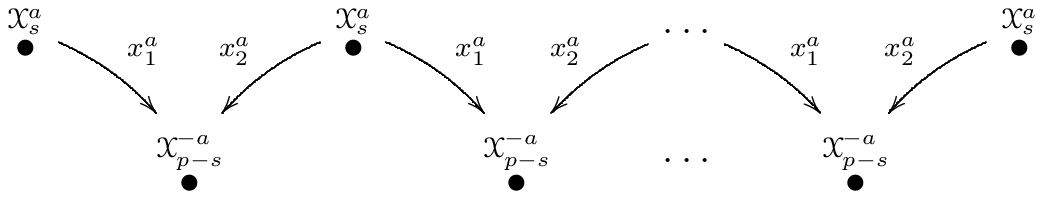}

  \noindent
  where $n$ is the number of $\repX^{a}_{s}$ modules.  A basis and the
  $\UresSL2$-action on $\modW^{a}_{s}(n)$ are explicitly described
  in~\bref{app:modW-base} in the example of $\modW^{a}_{s}(2)$.

\item[$\boldsymbol{\modM^{\pm}_{s}(n)}$] For $1\leq s \leq p{-}1$,
 $a\,{=}\,\pm\,$, and integer $n\geq2$, the module
 $\modM^{a}_{s}(n)$ is defined as\\*
 \includegraphics[bb=1.6in 9.6in 8in 10.5in, clip]{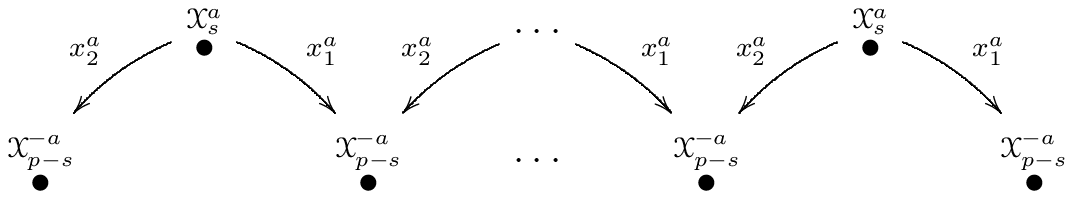}\\*
 where $n$ is the number of $\repX^{-a}_{p-s}$ modules.  A basis and
 the $\UresSL2$-action on $\modM^{a}_{s}(n)$ are explicitly described
 in~\bref{app:modM-base} in the example of $\modM^{a}_{s}(2)$.

\item[$\boldsymbol{\modO^{\pm}_{s}(n,z)}$] For $1\leq s \leq p-1$,
 $a\,{=}\,\pm$, integer $n\geq1$, and $z\in\oC\oP^1$, the module
 $\modO^{a}_{s}(n,z)$ is defined as
 \begin{center}
   \includegraphics[bb=1.6in 9.15in 8in 10.45in, clip]{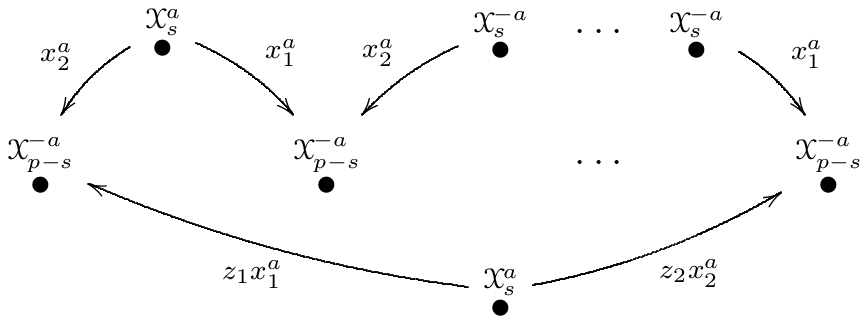}
 \end{center}
 where $n$ is the total number of the $\repX^{a}_{s}$ modules and
 $z=z_1:z_2\in\oC\oP^1$.  A basis and the $\UresSL2$-action on
 $\modO^{a}_{s}(n,z)$ are explicitly described in~\bref{app:modN-base}
 in the example of $\modO^{a}_{s}(1,z)$.
 
\item[$\boldsymbol{\modPr^{\pm}_{s}}$] For $1\leq s \leq p{-}1$ and
  $a\,{=}\,\pm\,$, the projective module $\modPr^{a}_{s}$ is defined
  in~\bref{app:proj-mod-base}; it has the following structure of
  subquotients (see~\bref{sec:proj-mod} below):
 \begin{equation}\label{schem-proj}
   \xymatrix@=12pt{
     &&\stackrel{\repX^a_{s}}{\bullet}
     \ar@/^/[dl]_{x^a_1} \ar@/_/[dr]^{x^a_2}&\\
     &\stackrel{\repX^{-a}_{p-s}}{\bullet}\ar@/^/[dr]_{x^a_2}&
     &\stackrel{\repX^{-a}_{p-s}}{\bullet}\ar@/_/[dl]^{x^a_1}\\
     &&\stackrel{\repX^a_{s}}{\bullet}&
   }
 \end{equation}
\end{description}
The following theorem shows that the modules listed above exhaust the
finite-dimensional indecomposable $\UresSL2$-modules.

\begin{Thmm}\label{thm:main}\mbox{}
  
  \begin{enumerate}
  \item The category $\lc$ has the decomposition
    \begin{equation*}
      \lc=\bigoplus_{s=0}^{p}\lc(s),
    \end{equation*}
    where each $\lc(s)$ is a full subcategory.
    
  \item Each of the subcategories $\lc(0)$ and~$\lc(p)$ is semisimple
    and contains precisely one irreducible module, $\repX^{+}_{p}$
    and~$\repX^{-}_{p}$ respectively.
    
  \item Each subcategory $\lc(s)$, $1\leq s \leq p{-}1$, contains
    precisely two irreducible modules~$\repX^{+}_{s}$
    and~$\repX^{-}_{p-s}$ and the following set of indecomposable
    modules:%%%\enlargethispage{12pt}
    \begin{itemize}
    \item the projective modules $\modPr^{+}_{s}$ and
    $\modPr^{-}_{p-s}$;
      
    \item three series of indecomposable modules:
      
      \begin{itemize}
      \item the modules $\modW^{+}_{s}(n)$ and $\modW^{-}_{p-s}(n)$
        for all integer $n\geq2$,
        
      \item the modules $\modM^{+}_{s}(n)$ and $\modM^{-}_{p-s}(n)$
        for all integer $n\geq2$,
        
      \item the modules $\modO^{+}_{s}(n,z)$ and
        $\modO^{-}_{p-s}(n,z)$ for all $z\in\oC\oP^1$ and integer
        $n\geq1$.
      \end{itemize}
    \end{itemize}
    This exhausts the list of indecomposable modules in $\lc(s)$.
  \end{enumerate}
\end{Thmm}

Next, for each subcategory $\lc(s)$ with $1\leq s\leq p-1$, the vector
space
\begin{equation*}
  \EXT_s=\bigoplus_{n\geq0}\Extn(\repX^{+}_{s}
  \oplus\repX^{-}_{p-s},
  \repX^{+}_{s}\oplus\repX^{-}_{p-s})
\end{equation*}
is an associative algebra with respect to the Yoneda product.  In the
following theorem, we describe the algebraic structure of~$\EXT_s$.
\begin{Thm}\label{ext-alg}
  The algebra $\EXT_s$ is generated by $x^\pm_i$, $i,j=1,2$, with the
  defining relations
  \begin{equation*}
    x^+_i x^+_j=x^-_ix^-_j=x^+_1x^-_2+x^+_2x^-_1=x^-_1x^+_2+x^-_2x^+_1=0.
  \end{equation*}
\end{Thm}
We note that the derived category of representations of $\EXT_s$ is
equivalent to the derived category of $\lc(s)$ (and hence to the
derived category of a full subcategory of $\algW(p)$-representations).

\subsection*{Notation}
We set
\begin{equation*}
  \q=e^{\frac{i\pi}{p}},
\end{equation*}
with $p=2$ in Sec.~\ref{sec:equiv} and $p\in\oZ_{\geq2}$ elsewhere.

In what follows, $\oZ_+$ denotes the set $\oN\cup\{0\}$.

\bigskip

This paper is organized as follows.  The proof of the classification
result in~\bref{thm:main} is given in Sec.~\ref{sec:proof}.  In
Sec.~\ref{sec:ext-proof}, we prove~\bref{ext-alg}.  The equivalence of
the $W$-algebra and quantum-group representation categories is proved
in Sec.~\ref{sec:equiv} for $p=2$.

Explicit constructions of some modules are given in
Appendix~\ref{app:mod}.  The necessary facts about quivers are
recalled in Appendix~\ref{app:quivers}.

\section{Proof of the classification theorem}\label{sec:proof}
We prove~\bref{thm:main} by straightforward investigation of the
structure of the category~$\lc$.  We fix an integer $p\geq2$.  The
strategy of the proof is as follows.  In~\bref{sec:cat-decomp}, using
the Casimir operator, we decompose the category~$\lc$ into a direct
sum of full subcategories~$\lc(s)$.  In~\bref{lem:ext-irrep}, we then
calculate~$\Ext$ between irreducible representations in each~$\lc(s)$.
In~\bref{sec:proj-mod}, we construct projective modules in each full
subcategory~$\lc(s)$.  We next note that all projective modules
in~$\lc$ are injective modules.  This information suffices to ensure
that indecomposable modules with semisimple length~$3$ are projective
modules and that there are no modules with semisimple length~$4$ or
more.  Therefore, to complete the proof of~\bref{thm:main}, it remains
to classify indecomposable modules with semisimple length~$2$.  We do
this in~\bref{sec:semisimple-2} using a correspondence between modules
with semisimple length $2$ and indecomposable representations of the
Kronecker quiver.

\subsection{Some conventions and definitions} In setting the notation
and recalling the basic facts about
$\overline{\mathscr{U}}_{\q}\equiv\UresSL2$ needed below, we largely
follow~\cite{[FGST]}.

\subsubsection{The restricted quantum group
  $\smash{\UresSL2}$}\label{subsec:qgroup} The Hopf algebra $\UresSL2$
is generated by $E$, $F$, and $K$ with the relations given
in~\bref{sec:Ures}.  In all $\UresSL2$-modules, we choose a basis such
that the generator $K$ acts diagonally (the action of $K$ in any
module in $\lc$ is diagonalizable because $K^{2p}{=}\one$).

Let $\cZ$ denote the center of $\UresSL2$ (it is
$(3p-1)$-dimensional~\cite{[FGST]}).

\subsubsection{Irreducible modules}\label{subsec:irrep} 
The irreducible $\UresSL2$-modules are labeled by their highest
weights $\q^{s-1}$, where $s\in\oZ/2p\oZ$.  We also parameterize the
same highest weights as $a\q^{s-1}$, where $a\,{=}\,\pm$ and $1\leq
s\leq p$.  Then for $1\leq s\leq p$, the irreducible module with the
highest weight $\pm \q^{s-1}$ is denoted by $\repX^{\pm}_{s}$.  The
dimension-$s$ module $\repX^{\pm}_{s}$ is spanned by elements
$\botpr_n^{\pm}$, $0\leq n\leq s{-}1$, where $\botpr_0^{\pm}$ is the
highest-weight vector and the algebra action is given~by
\begin{align*}%%%\label{basis-irrep}
  K \botpr_n^{\pm} &=
  \pm \q^{s - 1 - 2n} \botpr_n^{\pm},\\
  E \botpr_{n}^{\pm} &=
  \pm [n][s - n]\botpr_{n - 1}^{\pm},\\
  F \botpr_n^{\pm} &= \botpr_{n + 1}^{\pm},
\end{align*}
where we set $\botpr_{-1}^{\pm}=\botpr_{s}^{\pm}=0$ and use the
standard notation
\begin{equation*}
  [n] = \ffrac{\q^n-\q^{-n}}{\q-\q^{-1}}.
\end{equation*}
$\repX^{+}_{1}$ is the trivial module.

\subsubsection{Extensions}
Let $A$ and $C$ be left $\UresSL2$-modules.  We say that a short exact
sequence of $\UresSL2$-modules $0\to A\to B\to C\to 0$, or
equivalently,
\begin{equation*}
  A\injection B\surjection C,
\end{equation*}
is an \textit{extension} of $C$ by $A$, and we let $\Ext(C,A)$ denote
the set of equivalence classes (see, e.g.,~\cite{[M]}) of extensions
of $C$ by~$A$.

\subsubsection{Semisimple length of a module} Let $\modN$ be a
$\UresSL2$-module.  We define a \textit{semisimple filtration} of
$\modN$ as a tower of submodules
\begin{equation}\label{semi-filtr}
  \modN=\modN_0\supset\modN_1\supset\ldots\supset\modN_l=0
\end{equation}
such that each quotient $\modN_i/\modN_{i+1}$ is semisimple.  The
number $l$ is called the \textit{length} of the filtration.  In the
set of semisimple filtrations of $\modN$, there exists a filtration
with the minimum length~$\ell$.  We call~$\ell$ the \textit{semisimple
  length} of $\modN$.

Evidently, an indecomposable module has the semisimple length not less
than~$2$.  Any semisimple module has semisimple length~$1$.

\subsection{Decomposition of the category $\lc$}\label{sec:cat-decomp}
The recall the Casimir element
\begin{equation*}%%%\label{eq:casimir}
  \cas=EF+\ffrac{\q^{-1}K+\q K^{-1}}{(\q-\q^{-1})^2}=
  FE+\ffrac{\q K+\q^{-1}K^{-1}}{(\q-\q^{-1})^2}\in\cZ,
\end{equation*}
which satisfies the minimal polynomial relation (e.g.,
see~\cite{[FGST]})
\begin{equation*}%%%\label{Cas-relation}
  \Psi_{2p}(\cas)=0,
\end{equation*}
where
\begin{equation*}
  \Psi_{2p}(x) =
  (x-\beta_0)\,(x-\beta_p)\prod_{j=1}^{p-1}(x-\beta_j)^2, \qquad
  \beta_j=\ffrac{\q^j+\q^{-j}}{(\q-\q^{-1})^2},
\end{equation*}
and hence generates a $2p$-dimensional subalgebra in~$\cZ$.  This
relation yields a decomposition of the category $\lc$ into the direct
sum of full subcategories $\lc(j)$ such that $(\cas-\beta_j)$ acts
nilpotently on objects in $\lc(j)$.  Because $\beta_j\neq \beta_{j'}$
for $0\leq j\neq j'\leq p$, there are~$p\,{+}\,1$ full
subcategories~$\lc(s)$ for $0\,{\leq}\, s\,{\leq}\, p$.  This proves
part~(1) of Theorem~\bref{thm:main}.

Each subcategory $\lc(s)$ for $1 \leq s \leq p{-}1$ contains precisely
two irreducible modules $\repX^{+}_{s}$ and $\repX^{-}_{p-s}$, because
the Casimir element acts by multiplication with $\beta_s$ on precisely
these two irreducible modules.  The irreducible modules
$\repX^{+}_{p}$ and $\repX^{-}_{p}$ corresponding to the respective
eigenvalues $\beta_p$ and $\beta_0$ belong to the respective
categories $\lc(p)$ and $\lc(0)$.

To describe all isomorphism classes of indecomposable representations,
we first describe extensions between irreducible representations.
Inspection
shows that $\Ext(\repX^{a}_{s},$\linebreak[0]$\repX^{-a}_{p-s})
\cong\oC^2$ and that
all other $\Ext$ between irreducible representations in~$\lc(s)$
vanish.
We thus have the following lemma.

\begin{Lemma}\label{lem:ext-irrep}
  \addcontentsline{toc}{subsection}{\thesubsection.
    $\smash[t]{\protect\Ext{}}$ between irreducible representation}
  For $1 \leq s\leq p$ and $a,a'\,{=}\,\pm$, there are vector-space
  isomorphisms
  \begin{equation*}
    \Ext(\repX^{a}_{s},\repX^{a'}_{s'})\cong
    \begin{cases}
      \oC^2,\quad a'=-a, s'=p{-}s,\\
      0,\quad \text{otherwise}.
    \end{cases}
  \end{equation*}
\end{Lemma}
We also see that part~(2) and the statement in~\bref{thm:main} about
irreducible modules of $\lc(s)$ hold.

\begin{rem}
  The $\oC^2$ space in this lemma is the fundamental representation of
  the $s\ell(2)$ Lie algebra whereby $\UresSL2$ is extended to
  Lusztig's quantum group ``with divided powers.''  This is in fact
  another manifestation of the equivalence of the $\UresSL2$ and
  $\algW(p)$ representation categories: the occurrence of the
  $s\ell(2)$ Lie algebra is well known in the $(1,p)$ conformal field
  theory models, where it is an outer automorphism of the
  \textit{triplet} $W$-algebra $\algW(p)$ and hence acts on some of
  the $\algW(p)$-representations, and on extensions of two irreducible
  representations in particular.
\end{rem}

\subsubsection{Verma modules and bases in $\oC^2$} For each $1\leq
s\leq p$ and $a\,{=}\,\pm$, let $\Verma^a_{s}$ denote the Verma module
with the highest weight $a \q^{s-1}$, and let $\CVerma^a_{s}$ be the
contragredient Verma module with the highest weight $a\q^{-s-1}$ (it
is contragredient to the Verma
module~$\Verma^{-a}_{p-s}$).\footnote{The Verma $\Verma^a_{s}$ and the
  contragredient Verma $\CVerma^a_{s}$ modules are particular cases of
  the respective $\modO^{a}_{s}(1,z)$ modules with
  $z=z_1\,{:}\,z_2=1\,{:}\,0$ and $z=0\,{:}\,1$.}  We note that
$\Verma^\pm_{p}=\CVerma^\pm_{p}=\repX^{\pm}_{p}$.

For each $1\leq s\leq$\linebreak[]0$p-1$, we choose the bases
$\{x^+_i\}$ and $\{x^-_i\}$, $i=1,2$, in the respective spaces
$\oC^2=\Ext(\repX^{+}_{s},$\linebreak[0]$\repX^{-}_{p\,{-}\,s})$ and
$\oC^2=\Ext(\repX^{-}_{p\,{-}\,s},\repX^{+}_{s})$ such that $x^+_1$
corresponds to the Verma module $\Verma^+_{s}$, $x^+_2$ corresponds to
the contragredient Verma module $\CVerma^+_{s}$, $x^-_1$ corresponds
to~$\Verma^-_{p-s}$, and $x^-_2$ corresponds to~$\CVerma^-_{p-s}$.

\begin{rem}\label{rem:gluing}
  In terms of the action of the $\UresSL2$ generators $E$ and $F$, the
  extension
  $\repX^{+}_{s}\stackrel{x^+_1}{\longrightarrow}\repX^{-}_{p-s}$,
  giving the Verma module, is such that $F$ sends the
  \textit{lowest}-weight vector in the subquotient $\repX^{+}_{s}$ to
  the \textit{highest}-weight vector in the submodule
  $\repX^{-}_{p-s}$; similarly, the extension determined by $x^+_2$,
  giving the contragredient Verma module, is such that $E$ sends the
  highest-weight vector in $\repX^{+}_{s}$ to the lowest-weight vector
  in $\repX^{-}_{p-s}$.  The extensions determined by $x^{-}_i$,
  $i=1,2$, are described similarly.
\end{rem}

\subsection{Projective modules}\label{sec:proj-mod}
We next describe projective modules in each subcategory $\lc(s)$ for
$0\,{\leq}\, s \,{\leq}\,p$.  To construct these modules, we use that
any indecomposable projective module can be obtained as a projective
cover of an irreducible module.
\begin{prop}\label{prop:proj}\mbox{}
  
  \begin{enumerate}
  \item Projective modules of the subcategories $\lc(0)$ and $ \lc(p)$
    are the respective irreducible modules $\repX^{+}_{p}$ and
    $\repX^{-}_{p}$.
    
  \item For each $s$, $1\,{\leq}\, s \,{\leq}\,p{-}1$, projective
    modules of the subcategory~$\lc(s)$ are $\modPr^{+}_{s}$ and
    $\modPr^{-}_{p-s}$ \textup{(}see~\eqref{schem-proj}\textup{)}.
  \end{enumerate}
\end{prop}
\begin{proof}
  Part~(1) is clear because the subcategories $\lc(0)$ and $ \lc(p)$
  are semisimple.  To prove part~(2), we construct a projective cover
  of each irreducible module.  First, we construct a nontrivial
  extension of $\repX^{+}_{s}$ by the maximal number of irreducible
  modules $\repX^{-}_{p-s}$, that is, by $\repX^{-}_{p-s}\otimes
  \Ext(\repX^{+}_{s},\repX^{-}_{p-s})$, which is
  \begin{equation*}
    \repX^{-}_{p-s}\oplus\repX^{-}_{p-s}\injection\modM^{+}_{s}(2)
    \surjection\repX^{+}_{s}.
  \end{equation*}
  The indecomposable module $\modM^{+}_{s}(2)$ is described explicitly
  in~\bref{app:modM-base}.
  
  Next, to find the projective cover of $\repX^{+}_{s}$, we extend
  each summand in the submodule
  $\repX^{-}_{p-s}\oplus\repX^{-}_{p-s}\subset\modM^{+}_{s}(2)$ by
  $\repX^{+}_{s}\otimes \Ext(\repX^{-}_{p-s},\repX^{+}_{s})$.  The
  compatibility with the $\UresSL2$-algebra relations (with
  $F^p\,{=}\,0$ and $E^{p}\,{=}\,0$ in particular) leads to an
  extension corresponding to the module $\modPr^{+}_{s}$, which is the
  projective cover of $\repX^{+}_{s}$.  This module is described
  explicitly in~\bref{app:proj-mod-base}.
  
  A similar procedure gives the projective cover $\modPr^{-}_{p-s}$ of
  $\repX^{-}_{p-s}$.  The modules $\modM^{-}_{p-s}(2)$ and
  $\modPr^{-}_{p-s}$ are described in~\bref{app:modM-base}
  and~\bref{app:proj-mod-base} respectively.
  
  Therefore, in each subcategory $\lc(s)$ with $1\,{\leq}\,
  s\,{\leq}\, p{-}1$, we have two inequivalent projective modules
  $\modPr^{+}_{s}$ and $\modPr^{-}_{p-s}$, which cover the respective
  irreducible modules.  To show that these projective covers exhaust
  all indecomposable projective modules in the category~$\lc$, we
  recall that the left regular representation $\repReg$ of $\UresSL2$
  is a direct sum over all projective modules with the multiplicity of
  each given by the dimension of its simple quotient.  The dimension
  of the simple quotient of $\modPr^\pm_{s}$ is $s$
  (see~\bref{app:proj-mod-base}) and $\dim\repX^\pm_{p}=p$, and we
  therefore have
  \begin{equation*}
    \mathsf{Reg}\supset\bigoplus_{s=1}^{p-1} s\modPr^{+}_{s}
    \oplus
    \bigoplus_{s=1}^{p-1} s\modPr^{-}_{s}
    \oplus p\repX^{+}_{p}\oplus p\repX^{-}_{p}.
  \end{equation*}
  But counting the dimension of the right-hand side gives
  $\sum_{s=1}^{p-1}s\cdot 2p+\sum_{s=1}^{p-1}s\cdot 2p+p\cdot p+p\cdot
  p=2p^3=\dim\UresSL2$, and therefore the above inclusion is an
  equality, and part~(2) follows.
\end{proof}

\begin{rem}
  It also follows from~\bref{prop:proj} that the subquotient structure
  of the projective modules is the one shown in~\eqref{schem-proj}.
\end{rem}

Next, the ``dual'' procedure of constructing the injective modules in
each subcategory~$\lc(s)$ shows that they in fact coincide with the
projective modules (injective modules are contragredient to projective
modules, but the module obtained by inverting all arrows
in~\eqref{schem-proj} is evidently isomorphic to the original one).
We therefore have the following proposition.

\begin{prop}\label{prop:proj-inj}
  The projective modules in $\lc$ are injective.
\end{prop}

\begin{prop}\mbox{}

  \begin{enumerate}
  \item There are no modules in $\lc$ with semisimple length greater
    than~$3$.
    
  \item The only indecomposable modules with semisimple length~$3$ are
    the projective modules $\modPr^\pm_{s}$ with~$1\leq s\leq p{-}1$.
  \end{enumerate}
\end{prop}
\begin{proof}
  For $1\leq s\leq p-1$, let $\modN$ be a module such that
  $\modN_0/\modN_{1}$ in~\eqref{semi-filtr} contains the simple
  summand~$\repX^\pm_{s}$.  This means that there is a nonvanishing
  mapping~$\modPr^\pm_{s}\to\modN$ that covers~$\repX^\pm_{s}$.  There
  are two possibilities for the mapping~$\modPr^\pm_{s}\to\modN$:
  (i)~it is an embedding or (ii)~it has a nonvanishing kernel.  In
  case~(i), we recall that the projective module is at the same time
  injective and is therefore a direct summand in any module into which
  it is embedded.  In case~(ii), the kernel necessarily contains the
  submodule~$\repX^\pm_{s}$ of~$\modPr^\pm_{s}$.  But this means that
  the subquotient~$\repX^\pm_{s}$ of~$\modN$ actually belongs to a
  direct summand in~$\modN$ with the semisimple length less than or
  equal to~$2$.
\end{proof}

\subsection{Modules with semisimple length $2$}\label{sec:semisimple-2}
To complete the proof of part (3) in~\bref{thm:main}, it remains to
classify finite-dimensional modules with the semisimple length
$\ell\,{=}\,2$.

\subsubsection{The category $\lc^{(2)}(s)$} For $1 \leq s \leq p{-}1$,
let $\lc^{(2)}(s)$ be the full subcategory of $\lc(s)$ consisting of
$\UresSL2$-modules with semisimple length $\ell\leq2$.  Obviously, any
module in $\lc^{(2)}(s)$ can be obtained either by
\begin{itemize}
\item the extension of a direct sum of the irreducible modules
  $\repX^{+}_{s}$ by a direct sum of the irreducible modules
  $\repX^{-}_{p-s}$ via a direct sum of
  $x^+\!\in\!\Ext(\repX^{+}_{s},\repX^{-}_{p-s})$
\end{itemize}
or by
\begin{itemize}
\item the extension of a direct sum of the irreducible modules
  $\repX^{-}_{p-s}$ by a direct sum of the irreducible modules
  $\repX^{+}_{s}$ via a direct sum of
  $x^-\,{\in}\,\Ext(\repX^{-}_{p-s},\repX^{+}_{s})$.
\end{itemize}

For $\morX^+\,{\in}\,\Ext\left(\bigoplus_{j=1}^m\repX^{+}_{s},
  \bigoplus_{i=1}^n\repX^{-}_{p-s}\right)$, let
$\modI^{+}_{\morX^{+}}(m,n)\in\ob(\lc^{(2)}(s))$ denote the module
defined by the extension
\begin{equation}\label{general-gluing}
  \bigoplus_{j=1}^m\repX^{+}_{s}
  \xrightarrow{\morX^{+}}\bigoplus_{i=1}^n\repX^{-}_{p-s},
\end{equation}
For $ \morX^-\,{\in}\,\Ext\left(\bigoplus_{j=1}^m\repX^{-}_{p-s},
  \bigoplus_{i=1}^n\repX^{+}_{s}\right)$, we define the modules
$\modI^{-}_{\morX^{-}}(m,n)$ similarly.  We also set
$\modI^{\pm}_{0}(m,0)\!=\!\bigoplus_{j=1}^m\repX^{+}_{s}$ and
$\modI^{\pm}_{0}(0,n)\!=\!\bigoplus_{j=1}^n\repX^{-}_{p-s}$.\footnote{The
  modules $\modW^{+}_s(n), \modM^{+}_s(n)$, and $\modO^+_{s}(n,z)$
  in~\textbf{1.\theindsec} are particular cases of the
  $\modI^{+}_{\morX^{+}}(m,n)$ modules.}

We define two full subcategories $\lc^{(2),+}(s)$ and $\lc^{(2),-}(s)$
of the category $\lc^{(2)}(s)$ as follows.  An object of
$\lc^{(2),+}(s)$ is either a semisimple module or a module $\modN$
such that $\modN/\modN_1=\bigoplus_{j=1}^m\repX^{+}_{s}$ for some
$m\,{\in}\,\oN$, where $\modN_1$ is the maximal semisimple submodule
(the socle); in other words, an object of $\lc^{(2),+}(s)$ is the
$\modI^{+}_{\morX^{+}}(m,n)$ module for some $m,n\in\oZ_+$.  Objects
of $\lc^{(2),-}(s)$ are defined similarly with
$\modN/\modN_1=\bigoplus_{j=1}^m\repX^{-}_{p-s}$.  We note that
$\ob(\lc^{(2)}(s))\,{=}\,\ob(\lc^{(2),+}(s))\,
{\cup}\,\ob(\lc^{(2),-}(s))$.

\bigskip

We now reduce the classification of modules with semisimple length~$2$
to the classification of indecomposable representations of the
Kronecker quiver $\qK$.  The reader is referred to~\cite{[CB],[ARS]}
and Appendix~\bref{app:quivers} for the necessary facts about quivers.

\begin{lemma}\label{lemma:mod-Kronecker}
  Each of the categories $\lc^{(2),+}(s)$ and $\lc^{(2),-}(s)$ is
  equivalent to the category $\Rep(\qK)$ of representations of the
  Kronecker quiver $\qK$.
\end{lemma}
\begin{proof}
  The lemma is almost obvious if we note that two maps $\varepsilon$
  and $\bar\varepsilon$ in
  \begin{equation*}
    \xymatrix{
      {\repX^-_{p-s}}\ar@/_/[r]_{\varepsilon}
      \ar@/^/[r]^{\bar\varepsilon} &*{\;\modM^+_{s}(2)}
    }
  \end{equation*}
  make up a quiver in the category of $\UresSL2$-representations;
  taking the functors of $\justHom$ to each of the two categories then
  establishes the equivalence.  In somewhat more formal terms, the
  equivalence, e.g., between the categories $\lc^{(2),+}(s)$ and
  $\Rep(\qK)$ is given by the functor $\FunF$ that acts on objects as
  \begin{equation}\label{FunF-act}
    \FunF(\modI^+_{\morX^+}(m,n))=((V_0,V_1),(r_{01},\bar{r}_{01})),
  \end{equation}
  where
  \begin{align*}
    V_0&=\Hom(\modM^+_s(2),
    \modI^{+}_{\morX^{+}}(m,n))\,{=}\,\oC^m,\\
    V_1&=\Hom(\repPi_{p-s},\modI^{+}_{\morX^{+}}(m,n))\,{=}\,\oC^n
  \end{align*}
  and, for two linearly independent homomorphisms
  $\varepsilon,\bar{\varepsilon}
  \in\Hom(\repPi_{p-s},\modM^+_s(2))\,{=}\,\oC^2$, the two linear maps
  $r_{01},\bar{r}_{01}\in\HomC(V_0,V_1)$ are defined as
  \begin{equation*}
    r_{01}(\varphi)=\varphi\circ\varepsilon,\quad
    \bar{r}_{01}(\varphi)=\varphi\circ\bar{\varepsilon}
  \end{equation*}
  for each $\varphi\in V_0$, and with the natural action on morphisms.
  
  The existence of a functor $\FunG$ such that both $\FunG\FunF$ and
  $\FunF\FunG$ are the identity functors is evident from the
  definitions of the categories $\lc^{(2),+}(s)$ and $\Rep(\Kron)$.  
\end{proof}

Propositions~\bref{lemma:mod-Kronecker} and~\bref{prop:repr-Kron}
immediately imply the desired classification of finite-dimensional
$\UresSL2$-modules with semisimple length $\ell=2$, thus completing
the proof of~\bref{thm:main}.

\section{Proof of Theorem~\bref{ext-alg}}\label{sec:ext-proof}
The strategy of the proof is as follows.  In~\bref{sec:resolutions},
we first find the projective and injective resolutions of irreducible
representations.  In~\bref{lem:extn-irrep}, we then calculate~$\Extn$
between irreducible representations in each~$\lc(s)$ and
in~\bref{sec:ext-alg}, study the algebra~$\EXT$.

\subsection{Projective and injective resolution of irreducible
  modules}\label{sec:resolutions} Inspection shows that the mappings
defined in the following lemma give rise to a resolution.
\begin{lemma}\label{lem:proj-res-irrep}
  For each $1\,{\leq}\, s \,{\leq}\,p{-}1$, the module $\repX^{a}_{s}$
  has the projective resolution
  \begin{equation}\label{proj-res-irrep}
    \ldots\xrightarrow{\partial_3}\modPr^{a}_{s}
    \oplus\modPr^{a}_{s}\oplus\modPr^{a}_{s}
    \xrightarrow{\partial_2}\modPr^{-a}_{p-s}\oplus
    \modPr^{-a}_{p-s}\xrightarrow{\partial_1}\modPr^{a}_{s}
    \stackrel{\partial_0}{\surjection}\repX^{a}_{s}
  \end{equation}
  where for even $n$, the $n$th term is given by
  \begin{equation*}
    \xrightarrow{\partial_n}
    \underbrace{\modPr^{-a}_{p-s}\oplus\dots\oplus
      \modPr^{-a}_{p-s}}_{n}\xrightarrow{\partial_{n-1}}
  \end{equation*}
  with the boundary morphism given by the throughout mapping in
  \begin{equation*}
    \partial_n:\underbrace{\modPr^{a}_{s}
      \oplus\dots\oplus\modPr^{a}_{s}}_{n+1}
    \surjection\modW^{a}_{s}(n+1)\embedding
    \underbrace{\modPr^{-a}_{p-s}
      \oplus\dots\oplus\modPr^{-a}_{p-s}}_{n},
  \end{equation*}
  and for odd $n$, the $n$th term and the boundary morphism
  $\partial_n$ are given by changing $a$ to $-a$ and $s$ to $p{-}s$.
\end{lemma}

The statement dual to~\bref{lem:proj-res-irrep} is as follows.
\begin{lemma}
  For each $s$, $1\,{\leq}\, s \,{\leq}\,p{-}1$, the module
  $\repX^{a}_{s}$ has the injective resolution
  \begin{equation}\label{inj-res-irrep}
    \repX^{a}_{s}\stackrel{\delta_0}{\injection}
    \modPr^{a}_{s}
    \xrightarrow{\delta_1}\modPr^{-a}_{p-s}
    \oplus\modPr^{-a}_{p-s}\xrightarrow{\delta_2}\modPr^{a}_{s}
    \oplus\modPr^{a}_{s}\oplus\modPr^{a}_{s}\xrightarrow{\delta_3}\ldots
  \end{equation}
  where for even $n$, the $n$th term is given by
  \begin{equation*}
    \xrightarrow{\delta_{n-1}}
    \underbrace{\modPr^{-a}_{p-s}\oplus\dots\oplus
      \modPr^{-a}_{p-s}}_{n}\xrightarrow{\delta_n}
  \end{equation*}
  with the coboundary morphism given by the throughout mapping in
  \begin{equation*}
    \delta_n: \underbrace{\modPr^{-a}_{p-s}
    \oplus\dots\oplus\modPr^{-a}_{p-s}}_{n}\surjection
    \modM^{-a}_{p-s}(n+1)\embedding\underbrace{\modPr^{a}_{s}
    \oplus\dots\oplus\modPr^{a}_{s}}_{n+1},
  \end{equation*}
  and for odd $n$, the $n$th term and the coboundary morphism
  $\delta_n$ are given by changing $a$ to $-a$ and $s$ to $p{-}s$.
\end{lemma}

We now calculate $\Extn$ between irreducible modules.
\begin{Lemma}\label{lem:extn-irrep}
  \addcontentsline{toc}{subsection}{\thesubsection. \ $n$-Extensions
    between irreducible representations} For $1 \leq s \leq p-1$ and
  $a\,{=}\,\pm$, there are the vector-space isomorphisms
  \begin{align*}
    \Extn(\repX^{a}_{s},\repX^{a}_{s})&\cong
    \begin{cases}
      \oC^{n+1}, &n\;\text{is even},\\
      0,&n\;\text{is odd},\\
    \end{cases}
    \\
    \Extn(\repX^{a}_{s},\repX^{-a}_{p-s})&\cong
    \begin{cases}
      0,& n\;\text{is even},\\
      \oC^{n+1},& n\;\text{is odd}.
    \end{cases}
  \end{align*}
\end{Lemma}
\begin{proof}  
  Applying the contravariant functor $\Hom(-,\repX^{a}_{s})$ to the
  projective resolution of $\repX^{a}_{s}$ gives the cochain
  complex
  \begin{multline*}
    \bP:
    0\xrightarrow{\delta_0}\Hom(\modPr^a_{s},\repX^{a}_{s})
    \xrightarrow{\delta_1}
    \Hom(\modPr^{-a}_{p-s}\oplus\modPr^{-a}_{p-s},\repX^{a}_{s})
    \xrightarrow{\delta_2}\\
    \xrightarrow{\delta_2}\Hom(\modPr^a_{s}\oplus\modPr^a_{s}
    \oplus\modPr^a_{s},\repX^{a}_{s})
    \xrightarrow{\delta_3}\dots,
  \end{multline*}
  where for even $n$, the $(n\!+\!1)$th term is given by
  \begin{equation*}
    \xrightarrow{\delta_n}\Hom(
    \underbrace{\modPr^{a}_{s}\oplus\dots\oplus
      \modPr^{a}_{s}}_{n+1},\repX^{a}_{s})\cong\oC^{n+1}
    \xrightarrow{\delta_{n+1}}
  \end{equation*}
  and for odd $n$, the $(n\!+\!1)$th term is given by
  \begin{equation*}
    \xrightarrow{\delta_n}\Hom(
    \underbrace{\modPr^{-a}_{p-s}\oplus\dots\oplus
      \modPr^{-a}_{p-s}}_{n+1},\repX^{a}_{s})\cong0
    \xrightarrow{\delta_{n+1}}
  \end{equation*}
  It follows that coboundary morphisms $\delta_i=0$, $i\geq1$.  We
  then have
  \begin{equation*}
    \Extn(\repX^{a}_{s},\repX^{a}_{s})=
    H^n(\bP)=\ker(\delta_{n+1})/\im(\delta_{n})=
    \begin{cases}
      \oC^{n+1},& n\;\text{is even},\\
      0,& n\;\text{is odd}.
    \end{cases}
  \end{equation*}
  
  To calculate $\Extn(\repX^{a}_{s},\repX^{-a}_{p-s})$, we apply the
  functor $\Hom(-,\repX^{-a}_{p-s})$ to~\eqref{proj-res-irrep}.
\end{proof}

We next note that the Lusztig quantum group acts on the projective
resolution of $\repX^{a}_s$, and we have the following lemma.
\begin{lemma}\label{lemma:ext-sl2}
  The vector spaces $\Extn(\repX^{a}_{s},\repX^{a}_{s})$ and
  $\Extn(\repX^{a}_{s},\repX^{-a}_{p-s})$ are irreducible
  $s\ell(2)$-modules.
\end{lemma}

\subsection{The associative algebra $\justExt^{\bullet}$}
\label{sec:ext-alg}
Next, for each subcategory $\lc(s)$ with $1\leq s\leq p-1$, we
consider the total $\justExt$-group
\begin{equation*}
  \EXT\equiv\EXT_s
  =\bigoplus_{n\geq0}\Extn(\repX^{+}_{s}\oplus\repX^{-}_{p-s},
  \repX^{+}_{s}\oplus\repX^{-}_{p-s}).
\end{equation*}
{}By~\bref{lem:extn-irrep}, we have the vector-space isomorphism
\begin{equation*}
  \Extn(\repX^{+}_{s}\oplus\repX^{-}_{p-s},
  \repX^{+}_{s}\oplus\repX^{-}_{p-s})\cong\oC^{n+1}\oplus\oC^{n+1},
  \quad n\in\oZ_+.
\end{equation*}
We write $\EXT$ for $\EXT_s$ because its algebraic structure is
actually independent of $s$. The graded vector space $\EXT$ is
equipped with an associative algebra structure given by the Yoneda
product.

In what follows, abusing the notation, we just write a representative
of an element in $\Extn$ instead of the respective equivalence class.

We recall that a basis $\{x^+_i\}\cup\{x^-_i\}$, $i=1,2$, in
$\Ext(\repX^{+}_{s}\oplus\repX^{-}_{p-s},
\repX^{+}_{s}\oplus\repX^{-}_{p-s})\cong\oC^{2}\oplus\oC^{2}$ can be
chosen such that
\begin{alignat*}{2}
  x^{-}_1:\;
  \repX^{+}_{s}&\injection\Verma^-_{p-s}
  \surjection\repX^{-}_{p-s},&\quad
  x^{+}_1:\;\repX^{-}_{p-s}
  &\injection\Verma^+_{s}\surjection\repX^{+}_{s},\\
  x^{-}_2:\;\repX^{+}_{s}
  &\injection\CVerma^-_{p-s}\surjection\repX^{-}_{p-s},&\quad
  x^{+}_2:\;\repX^{-}_{p-s}
  &\injection\CVerma^+_{s}\surjection\repX^{+}_{s}.
\end{alignat*}

\begin{lemma}\label{lemma:Ext-alg-gen}
  The algebra $\EXT$ is generated by $x^\pm_i$, $i=1,2$.
\end{lemma}

\begin{proof}
  We first consider $\Extn(\repX^a_s,\repX^{a}_s)$ for even $n$.
  By~\bref{lemma:ext-sl2}, $\Extn(\repX^a_s,\repX^{a}_s)$ is an
  irreducible $s\ell(2)$-module.  The highest-weight vector $x_0$ of
  this module is given by
  \begin{equation*}
    x_0=\underbrace{x^{-a}_1 \dots x^a_1 x^{-a}_1 x^a_1}_{n}.
  \end{equation*}
  This cohomology class is nonzero because its restriction to the
  subalgebra $\oC[F]$ of $\UresSL2$ is nonzero, as can be directly
  verified.  Because of the $s\ell(2)$ action, it therefore follows
  that $\Extn(\repX^a_s,\repX^{a}_s)$ is generated by the $x^{\pm}_i$.
  For $\Extn(\repX^a_s,\repX^{-a}_{p-s})$ and odd $n$, the argument is
  similar.
\end{proof}

Taking the Yoneda product of the above short exact sequences, we
obtain the following exact sequences in the space
$\Yext{2}(\repX^{+}_{s},\repX^{+}_{s})\cong\oC^3$:
\begin{gather*}
  x^{-}_1x^{+}_1:\;\repX^{+}_{s}
  \injection\Verma^-_{p-s}\to\Verma^+_{s}\surjection\repX^{+}_{s},\\
  x^{-}_2x^{+}_2:\;\repX^{+}_{s}
  \injection\CVerma^-_{p-s}\to\CVerma^+_{s}\surjection\repX^{+}_{s},\\
  x^{-}_2x^{+}_1:\;\repX^{+}_{s}
  \injection\CVerma^-_{p-s}\to\Verma^+_{s}\surjection\repX^{+}_{s},\\
  x^{-}_1x^{+}_2:\;\repX^{+}_{s}
  \injection\Verma^-_{p-s}\to\CVerma^+_{s}\surjection\repX^{+}_{s}
\end{gather*}
and the exact sequences $x^{+}_1x^{-}_1, x^{+}_2x^{-}_2,
x^{+}_2x^{-}_1, x^{+}_1x^{-}_2$ in
$\Yext{2}(\repX^{-}_{p-s},\repX^{-}_{p-s})\cong\oC^3$.

\begin{lemma}\label{lemma:Ext-alg-rel}
  $x^-_1x^+_2+x^-_2x^+_1=x^+_1x^-_2+x^+_2x^-_1=0$.
\end{lemma}

\begin{proof}
  These relations follow from the existence of the projective modules.
  Explicitly, the argument for $x^-_1x^+_2+x^-_2x^+_1$ can be
  expounded as follows.  We consider the short exact sequence
  \begin{equation}\label{short-exact-seq-1}
    \repX^{-}_{p-s}\oplus\repX^{-}_{p-s}
    \injection\modM^{+}_{s}(2)\surjection\repX^{+}_{s},
  \end{equation}
  which gives rise to the long exact $\justExt$-sequence
  \begin{multline*}%%%\label{long-exact-seq-1}
    %%\xymatrix{
    0\to\Hom(\repX^{+}_{s},\repX^{+}_{s})
    \to\Hom(\modM^+_{s}(2),\repX^{+}_{s})\to
    \Hom(\repX^{-}_{p-s}\oplus\repX^{-}_{p-s},\repX^{+}_{s})
    \xrightarrow{\omega_0}\\
    \xrightarrow{\omega_0}\Ext(\repX^{+}_{s},\repX^{+}_{s})
    \to\Ext(\modM^+_{s}(2),\repX^{+}_{s})\xrightarrow{\varphi}
    \Ext(\repX^{-}_{p-s}\oplus\repX^{-}_{p-s},\repX^{+}_{s})
    \xrightarrow{\omega_1}\\
    \xrightarrow{\omega_1}\Yext{2}(\repX^{+}_{s},\repX^{+}_{s})\to\dots,
    %%}
  \end{multline*}
  where $\omega_0$ and $\omega_1$ are the connecting homomorphisms,
  given by the Yoneda product with the given short exact sequence
  \eqref{short-exact-seq-1} up to a sign.
  
  We then consider the short exact sequence
  $y^+\!\in\!\Ext(\modM^+_{s}(2),\repX^{+}_{s})$ given by
  \begin{equation*}
    \repX^{+}_{s}\injection\modPr^{+}_{s}\surjection\modM^{+}_{s}(2).
  \end{equation*}
  The extension $y^+$ of $\modM^{+}_{s}(2)$ defines the extension of
  the submodule
  $\repX^{-}_{p-s}\oplus\repX^{-}_{p-s}\subset\modM^+_{s}(2)$ as the
  exact sequence $\varphi(y^+)\in\Ext(\repX^{-}_{p-s}
  \oplus\repX^{-}_{p-s},\repX^{+}_{s})$ given by
  \begin{equation*}
    \repX^{+}_{s}\injection\modW^{-}_{p-s}(2)
    \surjection\repX^{-}_{p-s}\oplus\repX^{-}_{p-s}.
  \end{equation*}
  As a consequence of the identity $\omega_1\circ\varphi=0$, we have
  \begin{equation*}
    0=\omega_1\circ\varphi(y^+)=\repX^{+}_{s}
    \injection\modW^{-}_{p-s}(2)\to\modM^+_{s}(2)
    \surjection\repX^{+}_{s}
  \end{equation*}
  in $\Yext{2}(\repX^{+}_{s},\repX^{+}_{s})$ or, equivalently,
  \begin{equation*}
    0=\omega_1\circ\varphi(y^+)=x^{-}_1x^{+}_2 + x^{-}_2x^{+}_1.
  \end{equation*}
  That $x^+_1x^-_2+x^+_2x^-_1=0$ is established similarly with the
  help of the short exact sequence
  $\repX^{+}_{s}\oplus\repX^{+}_{s}\injection\modM^{-}_{p-s}(2)
  \surjection\repX^{-}_{p-s}$.
\end{proof}

{}From~\bref{lemma:Ext-alg-gen} and~\bref{lemma:Ext-alg-rel}, it
follows that the algebra $\EXT$ is generated by $x^{\pm}_i$, $i=1,2$,
and has the defining relations $x^{-}_1x^{+}_2 + x^{-}_2x^{+}_1=0$ and
$x^{+}_1x^{-}_2 + x^{+}_2x^{-}_1=0$.  We also set $x^+_i
x^+_j=x^-_ix^-_j=0$, $i,j=1,2$, because the Yoneda product is only
defined for two exact sequences if the end of one sequence coincides
with the beginning the other.  This completes the proof
of~\bref{ext-alg}.

\section{Equivalence of categories for $p=2$: an explicit
  construction}\label{sec:equiv} We now prove the equivalence of the
$\algW(p)$ and $\UresSL2$ representation categories for $p=2$.  The
proof considerably simplifies compared with the case of general $p$
because the $\algW(p)$-generators are explicitly constructed in terms
of symplectic fermions~\cite{[kaus3]}.  Extending the proof to
general~$p$ requires an appropriate generalization of the symplectic
fermions, achieved by introducing a multiplet of first-order conformal
fields carrying an $\UresSL2$-action.  The implication of the
equivalence is that the pictures in \textbf{1.\theindsec} equally well
describe $\algW(p)$-modules.

\subsection{$\algW(2)$ algebra and symplectic fermions}
For $p=2$, the $W$-algebra $\algW(2)$~\cite{[K-first],[GK2]}
(see~\textbf{1.\thewsec}) is conveniently described in terms of the
so-called symplectic fermions~\cite{[kaus3]} with
the OPE (in the holomorphic sector, which we consider in what follows)
\begin{equation*}
  \theta(z)\bar\theta(w)=\log(z-w).
\end{equation*}
The associated Virasoro algebra with central charge $c=-2$ is
generated by the energy-momentum tensor (with normal ordering
understood)
\begin{equation*}
  T(z)=-\dd\theta(z)\dd\bar\theta(z)
\end{equation*}
and the $\algW(2)$ algebra is generated by $T(z)$ and the two currents
given by
\begin{equation*}
  W^+(z)=\dd^2\theta(z)\dd\theta(z),
  \qquad
  W^-(z)=\dd^2\bar\theta(z)\dd\bar\theta(z).
\end{equation*}

The $(1,2)$ model contains two sectors of fields: generated by
$\theta(z)$ and $\bar\theta(z)$ from the operator $\one$ (periodic
boundary conditions) and from the so-called twist operator
$\mu_\pm(z)$ (antiperiodic boundary conditions); the subscript $\pm$
is for two \textit{quantum-group} components.  In the sector of
periodic boundary conditions, we have
\begin{equation}\label{periodic}
  \begin{split}
    \theta(z)&=\theta^*_0+\theta_0\log z
    -\sum_{n\neq0}\ffrac{1}{n}\,\theta_n z^{-n}\\
    \bar\theta(z)&=\bar\theta^*_0+\bar\theta_0\log z
    -\sum_{n\neq0}\ffrac{1}{n}\,\bar\theta_n z^{-n}.
  \end{split}
\end{equation}
In the sector of antiperiodic boundary conditions, we have
\begin{equation}\label{antiperiodic}
  \begin{split}
    \theta(z)&=\sum_{n\in\oZ+\half}\ffrac{1}{n}\,\theta_n z^{-n}\\
    \bar\theta(z)&=\sum_{n\in\oZ+\half}\ffrac{1}{n}\,\bar\theta_n z^{-n}.
  \end{split}
\end{equation}
The commutation relation between the modes in the sector of periodic
boundary conditions are (with $[{~},{~}]$ denoting supercommutator, in
this case anticommutator)
\begin{equation*}
  [\theta_n,\bar\theta_m]=\delta_{m+n,0},
  \quad
  [\theta_n,\theta_m]=0,
  \quad
  [\bar\theta_n,\bar\theta_m]=0,
\end{equation*}
and the only nonzero commutation relations involving $\theta^*_0$ and
$\bar\theta^*_0$ are\footnote{Indeed, because
  $\theta(ze^{2i\pi})=\theta(z)+2i\pi\theta_0$, it follows that the
  anticommutator $[2i\pi\theta_0,\bar\theta(w)]$ is given by the
  increment of the OPE $\theta(z)\bar\theta(w)$ as $z$ is continued
  via $z\to w + (z-w)e^{2i\pi}$, that is,
  $\log((z-w)e^{2i\pi})-\log(z-w)=2i\pi$.  This is independent of~$w$,
  and therefore $[\theta_0,\bar\theta^*_0]=1$,
  $[\theta_0,\bar\theta_n]=0$.}
\begin{equation*}
  [\theta_0,\bar\theta^*_0]=1,
  \quad
  [\bar\theta_0,\theta^*_0]=1.
\end{equation*}

\subsection{Bosonization}Another description of the model (as
in~\textbf{1.\thewsec}) is given in terms of a free bosonic field that
is a solution of the ``Riccati equation''
\begin{equation*}
  \half\dd\varphi\dd\varphi(z)+\half\dd^2\varphi(z)=T(z).
\end{equation*}
This equation has two solutions $\varphi^+(z)$ and $\varphi^-(z)$,
which are scalar fields with the OPEs
\begin{equation*}
  \varphi^\pm(z)\,\varphi^\pm(w)= \log(z{-}w)
\end{equation*}
and are related by
\begin{equation*}
  \varphi^+_0 = -\varphi^-_0,
\end{equation*}
where
\begin{equation*}
  \varphi^\pm_0 = \ffrac{1}{2i\pi}\oint dz\,\dd\varphi^\pm(z).
\end{equation*}
The symplectic fermions can then be described in either of the two
mutually nonlocal pictures, in terms of either $\varphi^+$ or
$\varphi^-$.  In the $\varphi^+$-picture, we have
\begin{equation}\label{plus-picture}
  \dd\theta(z)=e^{-\varphi^+(z)},
  \qquad\bar\theta(z)=e^{\varphi^+(z)}
\end{equation}
and in the $\varphi^-$-picture,
\begin{equation}
  \theta(z)=-e^{\varphi^-(z)},
  \qquad\dd\bar\theta(z)=e^{-\varphi^-(z)}.
\end{equation}

A special role is played by the zero modes
\begin{align*}%%%\label{theta0}
  \theta_0&=\ffrac{1}{2i\pi}\oint\!du\,\dd\theta(u)
  =\ffrac{1}{2i\pi}\oint\!du\,e^{-\varphi^+(u)}\\
  \intertext{and}
  \label{bar-theta0}
  \bar\theta_0&=\ffrac{1}{2i\pi}\oint\!du\,\dd\bar\theta(u)
  =\ffrac{1}{2i\pi}\oint\!du\,e^{-\varphi^-(u)}.
\end{align*}
In the $\varphi^+$-picture, $\theta_0$ is a screening operator and
$\bar\theta_0$ is a (nonlocal) contour-removal operator; in the
$\varphi^-$-picture, $\theta_0$ is a contour-removal operator and
$\bar\theta_0$ is a screening operator.

In the $\varphi^-$-picture, we recover the construction
in~\eqref{the-W} with
$W^-(z)=e^{-2\varphi^-(z)}=\dd^2\bar\theta(z)\dd\bar\theta(z)$ and
$S_+=\frac{1}{2i\pi}\oint du e^{2\varphi^-(u)}=\frac{1}{2i\pi}\oint du
\dd\theta(u)\theta(u)$.

Each of the bosonized pictures allows a component of the twist
operator to be expressed locally.  Starting with the
$\varphi^+$-picture for definiteness, we set
\begin{equation*}
  \mu_+(z)=e^{\frac{1}{2}\varphi^+(z)}.
\end{equation*}
The second component is then given by 
\begin{equation}\label{mu-pm}
  \mu_-(z)=\ffrac{1}{2i\pi}\int_{C_z} du e^{-\varphi^+(u)}
  e^{\frac{1}{2}\varphi^+(z)},
\end{equation}
where $C_z$ denotes a contour running from $-\infty$ to~$z$.  In the
$\varphi^-$-picture, we then have
\begin{equation*}
 \mu_-(z)\sim e^{\frac{1}{2}\varphi^-(z)}.
\end{equation*}
The two components are vectors in the two-dimensional representation
of the quantum group, see below.

\subsection{The quantum group} As usual, the screening and contour
removal operators generate a quantum group.  In the $(1,p)$
logarithmic conformal field theory models, it is just $\UresSL2$ at
$\q=e^{\frac{i\pi}{p}}$.

We have $\q=i$ for the $(1,2)$ model, with the quantum group
generators explicitly given by
\begin{equation}\label{def-EF}
  E=\theta_0,\qquad F=\pi\bar\theta_0 K,
\end{equation}
with $E^2=F^2=0$, where
\begin{equation}\label{def-K}
  K=e^{\pm i\pi\varphi^\pm_0}
  =e^{-%%%i\pi
    \half
    \oint\dd\theta\bar\theta}
  =e^{%%%i\pi
    \half
    \oint\theta\dd\bar\theta}.
\end{equation}
The zero modes $\theta^*_0$ and $\bar\theta^*_0$ can then be
considered functions on the quantum group.  We also have
\begin{equation}\label{commth}
  \theta_0\bar\theta_0+\bar\theta_0\theta_0=\ffrac{\one-K^2}{2i\pi}.
\end{equation}

\subsection{$\algW(2)$- and $\UresSL2$-representations}
We recall from~\cite{[kaus2],[GK2],[FHST]} that the algebra $\algW(2)$
has four irreducible representations, $\WrepX^\pm_{1}$ and
$\WrepX^\pm_{2}$.  The $\WrepX^\pm_{2}$ are Steinberg modules
(irreducible, Verma, and projective simultaneously).  We also need the
projective modules $\WrepP^{\pm}_{1}$ with the structure of
subquotients described as
\begin{equation*}%%%\label{schem-projW}
  \xymatrix@=12pt{
    &&\stackrel{\WrepX^{\pm}_{1}}{\bullet}
    \ar[dl]\ar[dr]&\\
    &\stackrel{\WrepX^{\mp}_{1}}{\bullet}\ar[dr]&
    &\stackrel{\WrepX^{\mp}_{1}}{\bullet}\ar[dl]\\
    &&\stackrel{\WrepX^{\pm}_{1}}{\bullet}&
  }
\end{equation*}  
The field corresponding to the highest-weight vector of $\WrepP^+_{1}$
is $\logid(z)=\theta(z)\bar\theta(z)$ (the ``logarithmic partner'' of
the operator~$\one$).

To construct the $\algW(2)$-representations in terms of the symplectic
fermions, we introduce the vacuum vector~$\vac$ corresponding to the
operator $\one$ (on which the symplectic fermions act with integer
modes, see~\eqref{periodic}); it is equivalently defined by the
highest-weight conditions
\begin{alignat*}{2}
  \theta_n\vac&=0,&\quad &n\in\oZ_{\geq0},\\
  \bar\theta_n\vac&=0,&\quad &n\in\oZ_{\geq0},\\
  K\vac&=\vac,
\end{alignat*}
where we simultaneously indicate the action of the $K$
generator~\eqref{def-K}.

Let $M$ be the module generated from $\vac$ by the integer modes
$\theta_{n\leq-1}$ and
$\bar\theta_{n\leq-1}$.  %%%and $\bar\theta_0^*$.  
The integer-moded module over the symplectic fermions is then
\begin{equation*}
  \bs{M}=M + \theta_0^* M
  + \bar\theta_0^* M + \bar\theta_0^*\theta_0^* M.
\end{equation*}

As a $\algW(2)$-module, $\bs{M}$ decomposes into a direct sum of two
projective modules, $\WrepP^+_{1}$, generated from
$\theta^*_0\bar\theta^*_0\vac$, and $\WrepP^-_{1}$, generated from
$\theta_{-1}\theta^*_0\bar\theta^*_0\vac$.

As a $\UresSL2$-module, $\bs{M}$ decomposes into an infinite sum of
projective modules.  For any state $\ket{v}=v\ket{0}\in M$, where $v$
is a homogeneous even-order polynomial in $\theta_{\leq-1}$ and
$\bar\theta_{\leq-1}$, the projective $\UresSL2$-module $\modPr^+_1$
(see~\eqref{schem-proj}) is spanned by
$\{\bar\theta^*_0\theta^*_0\ket{v}, \theta^*_0\ket{v},
-\bar\theta^*_0\ket{v},$\linebreak[0]$\ket{v}\}$.  Similarly, for any
state $\ket{v}=v\ket{0}\in M$, where $v$ is a homogeneous odd-order
polynomial in $\theta_{\leq-1}$ and $\bar\theta_{\leq-1}$, the
projective $\UresSL2$-module $\modPr^-_1$ is spanned by
$\{\bar\theta^*_0\theta^*_0\ket{v},$\linebreak[0]$
\theta^*_0\ket{v},$\linebreak[0]$-\bar\theta^*_0\ket{v}, \ket{v}\}$.

Next, let $\muvac\pm$ be the vector corresponding to the operator
$\mu_\pm$, on which the symplectic fermions act with half-integer
modes, see~\eqref{antiperiodic}.  It is equivalently defined by the
highest-weight conditions
\begin{equation*}
  \begin{alignedat}{2}
    \theta_n\muvac\pm&=0,& &n\in\half+\oZ_{\geq0},\\
    \bar\theta_n\muvac\pm&=0,& &n\in\half+\oZ_{\geq0},\\
    K\muvac\pm&=\pm i\muvac\pm.
  \end{alignedat}
\end{equation*}
Additionally, it satisfies the relations
\begin{equation*}
  \theta_0\muvac-=\bar\theta_0\muvac+=0.
\end{equation*}
The vectors $\muvac\pm$ are two basis vectors in the two-dimensional
representation of $\UresSL2$, with the screening operators~$\theta_0$
and~$\bar\theta_0$ acting as in~\eqref{mu-pm}, that is,
\begin{equation*}
  \theta_0\muvac+=\muvac-,
  \qquad
  \bar\theta_0\muvac-=\ffrac{i}{\pi}\muvac+.
\end{equation*}
We let $\bs{N}_\pm$ denote the module over the symplectic fermions
generated from $\muvac\pm$ by half-integer modes.

As a $\algW(2)$-module, each of $\bs{N}_\pm$ decomposes into two
Steinberg $\algW(2)$-modules $\WrepX^+_2$ and $\WrepX^-_2$.

As a $\UresSL2$-module, the direct sum of $\algW(2)$-modules
$\bs{N}_+\oplus\bs{N}_-$ decomposes into an infinite sum of Steinberg
modules.  The Steinberg module $\repX^+_2$ is spanned by
$\{v\muvac+,$\linebreak[0]$v\muvac-\}$ for any homogeneous even-order
polynomial $v$ in $\theta_{\leq-1}$ and $\bar\theta_{\leq-1}$;
similarly, the Steinberg module $\repX^-_2$ is spanned by
$\{v\muvac+,v\muvac-\}$ for any homogeneous odd-order polynomial $v$
in $\theta_{\leq-\half}$ and $\bar\theta_{\leq-\half}$.

\subsection{Morita equivalence and the proof of
  Theorem~\ref{Thm:equv}} Let $\Mmod$ denote the module
\begin{equation}\label{M-mod}
  \Mmod=\bs{M}\oplus\bs{N}_+\oplus\bs{N}_-.
\end{equation}
As a $\algW(2)$-module, it decomposes into a sum of \textit{all} the
projective modules $\WrepP^+_1$, $\WrepP^-_1$, $\WrepX^+_2$, and
$\WrepX^-_2$ with the multiplicities given by the dimensions of the
corresponding quantum-group modules,
\begin{align*}
  \Mmod&=\WrepP^+_1\boxtimes\repX^+_1\oplus
  \WrepP^-_1\boxtimes\repX^-_1
  \oplus\WrepX^+_2\boxtimes\repX^+_2\oplus
  \WrepX^-_2\boxtimes\repX^-_2.\\
  \intertext{As a $\UresSL2$-module, $\Mmod$ also decomposes into a
    sum of \textit{all} projective modules, with the
    ``multiplicities'' given by $\algW(2)$-modules:}
  \Mmod&=\WrepX^+_1\boxtimes\modPr^+_1\oplus
  \WrepX^-_1\boxtimes\modPr^-_1
  \oplus\WrepX^+_2\boxtimes\repX^+_2\oplus
  \WrepX^-_2\boxtimes\repX^-_2.
\end{align*}

The fact that $\Mmod$ is simultaneously a $\algW(2)$-module and a
$\UresSL2$-module and decomposes with respect to the projective
modules in each case allows establishing the equivalence of the
representation categories in the standard way.  We first construct a
functor $\bar{\mathcal{F}}:\rep{C}_2\to\mathfrak{C}_2$; it is actually
more convenient to first define the contravariant functor
\begin{equation*}
  \widetilde{\mathcal{F}}(\rep{Y})= \HomSL(\Mmod,\rep{Y}),
\end{equation*}
which gives a $\algW(2)$-module for any $\UresSL2$-module $\rep{Y}$,
and then compose it with the functor
$\boldsymbol{\mathsf{C}}:\mathfrak{C}_2\to\mathfrak{C}_2$ that
replaces each $\algW(2)$-module with its contragredient one.  Thus,
$\bar{\mathcal{F}} =\boldsymbol{\mathsf{C}}\widetilde{\mathcal{F}}:
\rep{C}_2\to\mathfrak{C}_2$ is a covariant functor.

We next construct the functor $\mathcal{F}:\mathfrak{C}_2\to\rep{C}_2$
similarly.  For any $\algW(2)$-module~$\WrepX$, we set
\begin{equation*}
  \mathcal{F}(\WrepX)
  =\mathcal{C}\,\HomW{2}(\Mmod,\WrepX)
\end{equation*}
(which is a $\UresSL2$-module), where
$\mathcal{C}:\rep{C}_2\to\rep{C}_2$ is the functor of taking the
contragredient module.

The functors $ \mathcal{F}$ and $ \bar{\mathcal{F}}$ send irreducible
representations into irreducible, and it therefore immediately follows
that
\begin{equation*}
  \mathcal{F}\bar{\mathcal{F}}\sim \mathrm{Id}_{\rep{C}_2},\quad
  \bar{\mathcal{F}}\mathcal{F}\sim \mathrm{Id}_{\mathfrak{C}_2}.
\end{equation*}

The natural braiding in conformal field theory, determined by contour
manipulations in the complex plane, coincides with the braiding
determined by the quantum-group $R$-matrix.  This is illustrated below
with several characteristic examples.  We first note that, as already
remarked in the introduction, the $\UresSL2$ quantum group at
$\q=e^{\frac{i\pi}{p}}$ is not quasitriangular, but is a subalgebra in
a larger quantum group generated by $E$, $F$, and $k$ such that
$k^2=K$, which \textit{is} quasitriangular~\cite{[FGST]}.  For $p=2$,
the universal $R$-matrix for this larger quantum group becomes
\begin{equation*}
  R=\ffrac{1}{8}\sum_{n,m=0}^7 \Bigl(
  i^{-\frac{nm}{2}}+2 i^{n-m-\frac{nm}{2}+1}E\tensor F
  \Bigr)k^n\tensor k^m,
\end{equation*}
where we choose $k=e^{+\half i\pi\varphi^+_0}$ to correspond to the
braiding matrix~$B^+$.

In view of~\eqref{def-EF}, the action of $R$ on vertex operators is
defined, which allows identifying $R$ with the~$B^+$ matrix.  Viewed
in this capacity, the $R$-matrix is rewritten as
\begin{equation*}
  R=\ffrac{1}{8}\sum_{n,m=0}^7 \Bigl(
  i^{-\frac{nm}{2}}+2\pi i^{n-m-\frac{nm}{2}+1}
  \theta_0\tensor \bar\theta_0\, e^{i\pi\varphi^+_0}
  \Bigr)e^{\half i\pi n\varphi^+_0}\tensor e^{\half i\pi m\varphi^+_0}.
\end{equation*}
This implies the braided relations
\begin{multline}\label{braiding}
  A(z)B(w)=RB(w)A(z)=\ffrac{1}{8}\sum_{n,m=0}^7
  i^{-\frac{nm}{2}}
  e^{\half i\pi n\varphi^+_0}B(w) e^{\half i\pi m\varphi^+_0}A(z)\\*
  -\ffrac{1}{8}\sum_{n,m=0}^7
  2\pi i^{-\frac{nm}{2}+1}e^{\half i\pi n\varphi^+_0}\theta_0B(w)
  e^{\half i\pi (m+2)\varphi^+_0}\bar\theta_0A(z)
\end{multline}
for any two fields $A(z)$ and $B(z)$, where $z$ and $w$ are real, with
$z>w$ (the relations are then analytically continued from the real
axis).  

In particular, we have
\begin{equation}\label{transposing}
  \theta(z)\bar\theta(w)=-\bar\theta(w)\theta(z)-2i\pi,\quad
  \bar\theta(z)\theta(w)=-\theta(w)\bar\theta(z),\qquad
  z>w.
\end{equation}
These relations can be illustrated by the procedure of permuting
operators in the complex plane.  In the first relation
in~\eqref{transposing}, in the $\varphi^+$-picture
(see~\eqref{plus-picture}), the operator $\theta(z)$ is nonlocal,
$\theta(z)=\int^z du e^{-\varphi^+(u)}$, which we represent with the
integration contour from infinity to its position,
\begin{equation*}
  \begin{picture}(210,60)(0,30)
    {\linethickness{1pt}
      \put(0,80){\line(1,0){200}}}
    \qbezier(200,80)(135,-10)(70,75)
    \put(70,75){\vector(-1,2){1}}
    \put(200,80){\circle*{2}} %% the point for \theta
    \put(205,80){$\theta(z)$}
    %% \put(135,0){\circle*{2}}
    %% \put(70,70){\circle*{2}}
    \put(110,75){\circle*{2}} %% the point for \bar\theta
    \put(115,68){$\bar\theta(w)$}
  \end{picture}
\end{equation*}
and taking it half the circumference (clockwise) around $w$ changes
the integration contour, resulting in the additional term~$-2i\pi$.
No such term occurs in the second relation in~\eqref{transposing},
represented as
\begin{equation*}
  \begin{picture}(210,60)(0,30)
    {\linethickness{1pt}
      \put(0,80){\line(1,0){110}}}
    \put(110,80){\circle*{2}} %% the point for \theta
    \put(115,68){$\theta(w)$}
    \qbezier(190,80)(135,-10)(70,75)
    \put(70,75){\vector(-1,2){1}}
    \put(190,80){\circle*{2}} %% the point for \bar\theta
    \put(195,80){$\bar\theta(z)$}
    %% \put(135,0){\circle*{2}}
    %% \put(70,70){\circle*{2}}
  \end{picture}
\end{equation*}

Similarly, for the field $\logid(z)=\theta(z)\bar\theta(z)$, which is
the logarithmic partner of $\one$, it follows from~\eqref{braiding}
that
\begin{equation*}
  \logid(z)\logid(w)=\logid(w)\logid(z)-2i\pi\theta(w)\bar\theta(z),
\end{equation*}
and for the twist operator, we have
\begin{gather*}
  \mu_\pm(z)\mu_\pm(w)=e^{\frac{i\pi}{4}}\mu_\pm(w)\mu_\pm(z),\quad
  \mu_+(z)\mu_-(w)=e^{\frac{i\pi}{4}}\mu_-(w)\mu_+(z),\\
  \mu_-(z)\mu_+(w)
  =e^{\frac{i\pi}{4}}(\mu_+(w)\mu_-(z)+2i\mu_-(w)\mu_+(z)).
\end{gather*}

The ribbon structure of $\UresSL2$ is given by the central element
\begin{equation*}
  \ribbon
  =\ffrac{1-i}{2\sqrt{2}}
  \Bigl((e^{\frac{i\pi}{4}}-2e^{-\frac{i\pi}{4}}KFE)(\one+K^2)
  +(K+2iFE)(\one-K^2)\Bigr).
\end{equation*}
Taking~\eqref{def-EF} into account, we obtain the element
\begin{multline*}
  \quad\ribbon
  =\ffrac{1-i}{2\sqrt{2}}\Bigl((e^{\frac{i\pi}{4}} +2\pi
  e^{i\pi(2\varphi^+_0-\frac{1}{4})}\bar\theta_0\theta_0)
  (\one+e^{2i\pi\varphi^+_0})\\*
  {}+2(\one-2i\pi\bar\theta_0\theta_0)e^{i\pi\varphi^+_0}
  (\one-e^{2i\pi\varphi^+_0})\Bigr),\quad
\end{multline*}
which acts on the space of states as
\begin{equation*}
  \ribbon=e^{2i\pi L_0}.
\end{equation*}

\section{Conclusions}
By the Kazhdan--Lusztig equivalence, the structure of indecomposable
modules pictured in \textbf{1.\theindsec} simultaneously describes the
structure of indecomposable $\algW(p)$-modules.  The projective
resolution in~\eqref{proj-res-irrep} is also applicable to the
$\algW(p)$ algebra, and may therefore be used in conformal field
theory as well as the known (Felder-type) resolutions are used in
rational conformal field theories.

As an application of the algebra $\EXT$ described in~\bref{ext-alg},
we note the following.  For a subcategory $\lc(s)$ in the category
$\lc$ (considered as either the $\algW(p)$- or
$\UresSL2$-representation category), the derived category of $\lc(s)$
is equivalent to the derived category of coherent sheaves on a
noncommutative extension of $\oC\oP^1$~\cite{[ABZ]}.  The algebra
$\EXT$ is then the coordinate ring of this noncommutative extension.
 
In addition to establishing the equivalence of categories, the
bimodule $\Mmod$ constructed in~\eqref{M-mod} has another application,
to the construction of the \textit{full} conformal field
theory~\cite{[MS],[fuRs4],[HK]}.  Given $\Mmod$ (a module over the
quantum group and the $\algW(p)$ algebra assigned to the holomorphic
sector) and its second copy $\bar\Mmod$ where the $\algW(p)$-algebra
action is replaced with its action in the antiholomorphic sector, the
standard construction of the full conformal field theory is to take
the tensor product $\Mmod\tensor_{\qgr}\bar\Mmod$ over the quantum
group $\qgr=\UresSL2$.  Other possibilities are to take
$\Mmod\tensor_{\qgr}\!N\!\tensor_{\qgr}\bar\Mmod$, where $N$ is any
$\qgr$-bimodule with an associative product $N\tensor N\to N$
consistent with the bimodule structure.  The bimodule $\Mmod$ may also
be used to obtain (the space dual to) the quantum group center as
``extended characters'' in logarithmic conformal field theory models.
The ``extended characters'' can be obtained by taking the trace of
$e^{2i\pi\tau L_0}$ over the bimodule $\Mmod$ (over its
vertex-operator-algebra component), which should give a version of the
general construction in~\cite{[My3]}.

\subsubsection*{Acknowledgments} We are grateful to S.~Ariki, T.~Miwa
and I.V.~Tyutin for useful remarks and to J.~Fuchs for valuable
discussions and comments on the literature.  Part of the paper was
written during our stay at Kyoto University, and we are grateful to
T.~Miwa for the kind hospitality extended to us.  This paper was
supported in part by the RFBR Grants~04-01-00303, 05-02-17451, and LSS
1578.2003.2 and the RFBR--JSPS Grant~05-01-02934YaF\_a.

\appendix

\section{$\UresSL2$-modules}\label{app:mod}

\subsection{The modules $\modW^{\pm}_{s}(2)$, $\modM^{\pm}_{s}(2)$,
  and $\modO^{\pm}_{s}(1,z)$}\label{app:modW-pm-base}

\subsubsection{The module $\modW^{{\pm}}_{s}(2)$}\label{app:modW-base}
Let $s$ be an integer $1\leq s\leq p{-}1$ and $a=\pm$.  The module
$\modW^{a}_{s}(2)$ (see~\textbf{1.\theindsec}) has three
subquotients: $\repX^{a}_{s}$ (``left''), $\repX^{-a}_{p-s}$
(``bottom''), and $\repX^{a}_{s}$ (``right'').  The basis can be
chosen as
\begin{equation*}
  \{\botpr_n,\toppr_n\}_{0\le n\le s-1}
  \cup\{\leftpr_k\}_{0\le k\le p-s-1},
\end{equation*}
where $\{\botpr_n\}_{0\le n\le s-1}$ is the basis corresponding to the
right $\repX^{a}_{s}$ module, $\{\toppr_n\}_{0\le n\le s-1}$ to the
left $\repX^{a}_{s}$, and $\{\leftpr_k\}_{0\le k\le p-s-1}$ to the
bottom $\repX^{-a}_{p-s}$ module, with the $\UresSL2$-action given by
\begin{alignat*}{3}
  K\botpr_n&=aq^{s-1-2n}\botpr_n,& \quad
  K\toppr_n&=aq^{s-1-2n}\toppr_n,& \quad &0\le n\le s{-}1,\\
  K\leftpr_k&=-aq^{p-s-1-2k}\leftpr_k,& \quad
  &0\le k\le p{-}s{-}1,\\
  E\leftpr_k&=-a[k][p-s-k]\leftpr_{k-1},& \quad 0\le
  k&\le p{-}s{-}1 \quad(\text{with}\quad\leftpr_{-1}\equiv0),
  \kern-60pt
\end{alignat*}
\begin{align*}
  E\botpr_n&=
  \begin{cases}
    a[n][s-n]\botpr_{n-1}, &1\le n\le s{-}1,\\
    \leftpr_{p-s-1}, & n=0,\\
  \end{cases}\\
  E\toppr_n&=a[n][s-n]\toppr_{n-1},
  \quad 0\le n\le s{-}1\quad(\text{with}\quad\toppr_{-1}\equiv0),\\
  \intertext{and} F\botpr_n&=\botpr_{n+1}, \quad 0\le
  n\le s{-}1
  \quad(\text{with}\quad\botpr_s\equiv0),\\
  F\toppr_n&=
  \begin{cases}
    \toppr_{n+1}, &0\le n\le s{-}2,\\
    \leftpr_0, & n=s{-}1,
  \end{cases}\\
  F\leftpr_k&=\leftpr_{k+1}, \quad 0\le k\le p{-}s{-}1
  \quad(\text{with}\quad\leftpr_{p-s}\equiv0).
\end{align*}

\subsubsection{The module $\modM^{{\pm}}_{s}(2)$}
\label{app:modM-base} Let $s$ be an integer $1\leq
s\leq p{-}1$ and $a=\pm$.  The module $\modM^{a}_{s}(2)$, which is the
extension
\begin{equation}\label{schem-M2}
  \xymatrix@=12pt{
    &\stackrel{\repX^{+}_{s}}{\bullet}\ar@/^/[dl]_{x^+_1}
    \ar@/_/[dr]^{x^+_2}&&\\
    \stackrel{\repX^{-}_{p-s}}{\bullet}&&
    \stackrel{\repX^{-}_{p-s}}{\bullet}&
  }
\end{equation}
has the basis
\begin{equation*}
  \{\botpr_n\}_{0\le n\le s-1}
  \cup\{\leftpr_k,\rightpr_k\}_{0\le k\le p-s-1},
\end{equation*}
where $\{\botpr_n\}_{0\le n\le s-1}$ is the basis
corresponding to the top module $\repX^{a}_{s}$ in~\eqref{schem-M2},\\
$\{\leftpr_k\}_{0\le k\le p-s-1}$ to the left $\repX^{-a}_{p-s}$, and
$\{\rightpr_k\}_{0\le k\le p-s-1}$ to the right $\repX^{-a}_{p-s}$,
with the $\UresSL2$-action given by
\begin{alignat*}{3}
  K\botpr_n&=aq^{s-1-2n}\botpr_n,& \quad &0\le n\le s{-}1,\\
  K\leftpr_k&=-aq^{p-s-1-2k}\leftpr_k,& \quad
  K\rightpr_k&=-aq^{p-s-1-2k}\rightpr_k, \quad &0\le k\le
  p{-}s{-}1,\\
  E\leftpr_k&=-a[k][p-s-k]\leftpr_{k-1},& \quad
  E\rightpr_k&=-a[k][p-s-k]\rightpr_{k-1},\quad &0\le k\le
  p{-}s{-}1\\*
  &&&(\text{with}\quad\leftpr_{-1}\equiv\rightpr_{-1}\equiv0),& 
  %% \kern-60pt
\end{alignat*}
\begin{align*}
  E\botpr_n&=
  \begin{cases}
    a[n][s-n]\botpr_{n-1}, &1\le n\le s{-}1,\\
    \leftpr_{p-s-1}, & n=0,\\
  \end{cases}\\
  \intertext{and}
  F\leftpr_k&=\leftpr_{k+1}, \quad
  F\rightpr_k=\rightpr_{k+1}, \quad 0\le k\le p{-}s{-}1
  &\quad(\text{with}\quad\leftpr_{p-s}\equiv\rightpr_{p-s}\equiv0),\\
  F\botpr_n&=
  \begin{cases}
    \botpr_{n+1}, &0\le n\le s{-}2,\\
    \rightpr_0, & n=s{-}1.
  \end{cases}
\end{align*}

\subsubsection{The module $\modO^{\pm}_{s}(1,z)$}\label{app:modN-base}
Let $s$ be an integer $1\leq s\leq p{-}1$, $a=\pm$, and
$z\,{=}\,\abp$.  For $n=1$, the diagram for the $\modO^{a}_{s}(1,z)$
module in~\textbf{1.\theindsec} %%%~\bref{sec:the-modules} 
takes the form
\begin{equation*}
  \xymatrix@=12pt{
    \stackrel{\repX^{a}_{s}}{\bullet}\ar@/^3ex/[dr]^{z_1x^a_1}
    \ar@/_3ex/[dr]_{z_2x^a_2}&\\ &\stackrel{\repX^{-a}_{p-s}}{\bullet}
  }
\end{equation*}
In terms of the respective bases $\{\botpr_n\}_{0\le n\le s-1}$ and
$\{\leftpr_k\}_{0\le k\le p-s-1}$ in $\repX^{a}_{s}$ and
$\repX^{-a}_{p-s}$, the $\modO^{a}_{s}(1,z)$ module has the basis
\begin{equation}\label{Nbasis}
  \{\botpr_n\}_{0\le n\le s-1}\cup\{\leftpr_k\}_{0\le k\le p-s-1}
\end{equation}
with the $\UresSL2$-action given by
\begin{align*}
  K\botpr_n&=aq^{s-1-2n}\botpr_n, \quad
  K\leftpr_k=-aq^{p-s-1-2k}\leftpr_k,
  \quad 0\le n\le s{-}1,~0\le k\le p{-}s{-}1,\\
  E\botpr_n&=
  \begin{cases}
    a[n][s-n]\botpr_{n-1}, &1\le n\le s{-}1,\\
    z_2 \leftpr_{p-s-1}, &n=0,\\
  \end{cases}
  \\
  E\leftpr_k&=-a[k][p-s-k]\leftpr_{k-1}, \quad 0\le k\le p{-}s{-}1,\\
  \intertext{where we set $\leftpr_{-1}=0$, and} F\botpr_n&=
  \begin{cases}
    \botpr_{n+1}, &0\le n\le s{-}2,\\
    z_1 \leftpr_0, &n=s{-}1,\\
  \end{cases}
  \\
  F\leftpr_k&=\leftpr_{k+1}, \quad 0\le k\le p{-}s{-}1,
\end{align*}
where we set $\leftpr_{p-s}=0$.  That this defines a $\UresSL2$-module
is verified by direct calculation.

\subsection{The projective modules
  $\modPr^{\pm}_{s}$}\label{app:proj-mod-base}
Let $s$ be an integer $1\leq s\leq p{-}1$ and $a=\pm$.  The projective
module $\modPr^a_{s}$ has the basis
\begin{equation*}
  \{\botpr_n,\toppr_n\}_{0\le n\le s-1}
  \cup\{\leftpr_k,\rightpr_k\}_{0\le k\le p-s-1},
\end{equation*}
where $\{\toppr_n\}_{0\le n\le s-1}$ is the basis corresponding to the
top module $\repX^{a}_{s}$ in~\eqref{schem-proj},\\ $\{\botpr_n\}_{0\le
  n\le s-1}$ to the bottom $\repX^{a}_{s}$, $\{\leftpr_k\}_{0\le k\le
  p-s-1}$ to the right $\repX^{-a}_{p-s}$, and $\{\rightpr_k\}_{0\le
  k\le p-s-1}$ to the left module $\repX^{-a}_{p-s}$, with the
$\UresSL2$-action given by
\begin{alignat*}{3}
  K\leftpr_k&=-aq^{p-s-1-2k}\leftpr_k,& \quad
  K\rightpr_k&=-aq^{p-s-1-2k}\rightpr_k,&
  \quad &0\le k\le p{-}s{-}1,\\
  K\botpr_n&=aq^{s-1-2n}\botpr_n,& \quad
  K\toppr_n&=aq^{s-1-2n}\toppr_n,& \quad &0\le n\le s{-}1,\\
  E\leftpr_k&=-a[k][p-s-k]\leftpr_{k-1},& \quad 0\le
  k&\le p{-}s{-}1 \quad(\text{with}\quad\leftpr_{-1}\equiv0),
  \kern-60pt
\end{alignat*}
\begin{align*}
  E\rightpr_k&=
  \begin{cases}
    -a[k][p-s-k]\rightpr_{k-1}, &1\le k\le p{-}s{-}1,\\
    \botpr_{s-1}, & k=0,\\
  \end{cases}
  \\
  E\botpr_n&=a[n][s-n]\botpr_{n-1},
  \quad 0\le n\le s{-}1\quad(\text{with}\quad\botpr_{-1}\equiv0),\\
  E\toppr_n&=
  \begin{cases}
    a[n][s-n]\toppr_{n-1}+\botpr_{n-1}, &1\le n\le s{-}1,\\
    \leftpr_{p-s-1}, & n=0,\\
  \end{cases}
  \\
  \intertext{and} F\leftpr_k&=
  \begin{cases}
    \leftpr_{k+1}, &0\le k\le p{-}s{-}2,\\
    \botpr_0, & k=p{-}s{-}1,\\
  \end{cases}
  \\
  F\rightpr_k&=\rightpr_{k+1}, \quad 0\le k\le
  p{-}s{-}1
  \quad(\text{with}\quad\rightpr_{p-s}\equiv0),\\
  F\botpr_n&=\botpr_{n+1}, \quad 0\le n\le s{-}1
  \quad(\text{with}\quad\botpr_s\equiv0),\\
  F\toppr_n&=
  \begin{cases}
    \toppr_{n+1}, &0\le n\le s{-}2,\\
    \rightpr_0, & n=s{-}1.
  \end{cases}
\end{align*}

\section{The Kronecker quiver}\label{app:quivers}
We here recall basic notions about quivers~\cite{[CB],[ARS]}.

\subsection{Quivers and their representations}\label{sec:quivers}
A quiver $Q$ is an oriented graph, that is, a quadruple $Q=(I, A, s,
t)$, consisting of a finite set $I$ of vertices, a finite set $A$ of
oriented edges (arrows), and two maps $s$ and $t$ from $A$ to $I$.  An
oriented edge $a \in A$ starts at the vertex $s(a)$ and terminates at
$t(a)$.

A \textit{representation} of a quiver $Q$ (over $\oC$) is a collection
of finite-dimensional vector spaces $V_i$ over $\oC$, one for each
vertex~$i \in I$, and $\oC$-linear maps $r_{ij}: V_i \to V_j$, one for
each oriented edge $\stackrel{i}{\bullet}\xrightarrow{a}
\stackrel{j}{\bullet}$.  \ The \textit{dimension} of a representation
$\rho$ of $Q$ is an element of $\oZ[I]$ given by the dimensions of
$V_i$, $i\in I$: $\dim(\rho)= \sum_{\substack{i\in
    I}}\dim_{\oC}(V_i)i$.

A \textit{morphism} from a representation $\rho$ of a quiver $Q$ to
another representation $\rho'$ of $Q$ is an $I$-graded $\oC$-linear
map $\phi=\bigoplus_{i\in I}\phi_i:\bigoplus_{i\in
  I}V_{i}\to\bigoplus_{i\in I}V_{i}'$ satisfying $r_{ij}'\phi_i=\phi_j
r_{ij}$ for each oriented edge $\stackrel{i}{\bullet}\xrightarrow{a}
\stackrel{j}{\bullet}$.  This gives the category of representations of
the quiver $Q$, to be denoted by $\Rep(Q)$ in what follows.

If a quiver $Q$ has no oriented cycles, isomorphism classes of simple
objects in $\Rep(Q)$ are in a one-to-one correspondence with vertices
of $Q$.  The simple object corresponding to a vertex $i\in I$ is given
by the vector spaces
\begin{align*}
  V_j&=
 \begin{cases}
   \oC,\quad j\!=\!i,\\
   0,\quad \text{otherwise}
 \end{cases}\\
 \intertext{and $\oC$-linear maps}
 r_{ij}&=0,\quad \text{for all}\; i,j\in I.
\end{align*}
A quiver is said to be of \textit{finite} type if the underlying
nonoriented graph is a Dynkin graph of finite type.  Similarly, a
quiver is said to be of \textit{affine} type if the underlying
nonoriented graph is a Dynkin graph of affine type.  A quiver is said
to be of \textit{simply laced} type if it does not have a pair of
vertices connected by more than one arrow.

\subsection{The category $\Rep(\qK)$ of representations of the
  Kronecker quiver}\label{subsec:rep-Kronecker} 
The Kronecker quiver $\qK$ has two vertices connected by two edges,
which are directed in the same way, that is, $\qK =
(\{0,1\},\{f,g\},s,t)$, where $s(f)=s(g)=0$ and $t(f)=t(g)=1$:
\begin{equation*}
  \xymatrix{
    {\stackrel{0}{\bullet}}\ar@/_/[r]_g
    \ar@/^/[r]^f &{\stackrel{1}{\bullet}}
  }
\end{equation*}

A representation $\rho$ of $\qK$ is a quadruple
$((V_0,V_1),(r_{01},\bar{r}_{01}))$ consisting of two $\oC$-linear
spaces $V_0$ and $V_1$ and two linear maps $r_{01},\bar{r}_{01}\in
\mathrm{Hom}_{\oC}(V_0,V_1)$.  The dimension of~$\rho$ is given by
$\dim(\rho)=(\dim_{\oC}(V_0), \dim_{\oC}(V_1))$.  Simple objects in
the category $\Rep(\Kron)$ are given by the two representations
$\rho_0=((\oC,0),(0,0))$ and $\bar{\rho}_0=((0,\oC),(0,0))$.  We now
recall the classification of indecomposable representations of the
Kronecker quiver $\qK$, summarized in~\bref{prop:repr-Kron} below.

There is a correspondence between indecomposable representations of a
quiver and the set $\Delta_{+}$ of positive roots of the Lie algebra
corresponding to the Dynkin graph associated with the quiver.  This
correspondence is one-to-one for a quiver of simply laced finite type
\cite{[Gab],[BGP]}.  Namely, a representation $\rho$ of a quiver $Q$
is indecomposable if and only if $\dim(\rho)\in\Delta_{+}$ and,
conversely, for every $\alpha\in\Delta_{+}$, there is, up to an
isomorphism, a unique indecomposable representation $\rho$ of the
quiver $Q$ such that $\dim(\rho)=\alpha$.

For quivers of affine type, the bijection between isomorphism classes
of indecomposable representations and positive roots of the
corresponding affine Lie algebra holds only for positive real roots
$\alpha\in\Delta^{\text{re}}_{+}$ \cite{[Naz],[DF]}.  For each positive
imaginary root $\alpha\in\Delta^{\text{im}}_{+}$ there exists an
uncountable set $\setI_{\alpha}$ of nonisomorphic indecomposable
representations corresponding to $\alpha$.  Actually, the set
$\setI_{\alpha}$ is independent of
$\alpha\in\Delta^{\text{im}}_{+}$~\cite{[DR]}.  If $\alpha\in\oZ[I]$ is
not in $\Delta_+$, then the set of indecomposable representations with
the dimension $\alpha$ is empty.

The nonoriented graph associated with the Kronecker quiver $\Kron$ is
the extended Dynkin graph $\Aone$.  It is well known that
$\alpha\in\Delta_+$ is a positive real root of $\Aone$ if
$\alpha=(n+1,n)$ or $\alpha=(n,n+1)$, and $\alpha\in\Delta_+$ is a
positive imaginary root if $\alpha=(n+1,n+1)$, where
$n\,{\in}\,\oZ_{+}$.  In particular, $\alpha_0=(1,0)$ and
$\bar{\alpha}_0=(0,1)$ are the simple roots and $\delta=(1,1)$ is the
first imaginary root.  The simple roots $\alpha_0$ and
$\bar{\alpha}_0$ correspond to the respective simple objects
$\rho_0=((\oC,0),(0,0))$ and $\bar{\rho}_0=((0,\oC),(0,0))$ in the
category $\Rep(\Kron)$.  The other real roots are in a one-to-one
correspondence with indecomposable representations
$\rho_n=((\oC^{n+1},\oC^n),(R,\bar{R}))$ and
$\bar{\rho}_n=((\oC^{n},\oC^{n+1}),(R^t,\bar{R}^t))$, where $R$ and
$\bar{R}$ are the \textit{rectangular} matrices
\begin{equation*}
  R=
  \begin{pmatrix}
    1 & 0 & \dots & 0 & 0 \\
    0 & 1 & \dots & 0 & 0 \\
    \hdotsfor{5} \\
    0 & 0 & \dots & 0 & 0 \\
    0 & 0 & \dots & 1 & 0 
  \end{pmatrix} , 
  \quad
  \bar{R}=
  \begin{pmatrix}
    0 & 1 & \dots & 0 & 0 \\
    0 & 0 & \dots & 0 & 0 \\
    \hdotsfor{5} \\
    0 & 0 & \dots & 1 & 0 \\
    0 & 0 & \dots & 0 & 1 
  \end{pmatrix}.
\end{equation*}
The first imaginary root $\delta$ corresponds to the family
$\setI_{\delta}$ of indecomposable representations
$((\oC,\oC),(r_{01},\bar{r}_{01}))$ for some
$r_{01},\bar{r}_{01}\in\HomC(\oC,\oC)$, parameterized by $\oC\oP^1$.
But the set $\setI_{\alpha}$ for a positive imaginary root~$\alpha$ is
independent of~$\alpha$, and therefore $\setI_{\alpha}=\oC\oP^1$ for
any imaginary root $\alpha$.  We summarize these results in the
following well-known proposition (see, e.g., \cite{[FMV]}), first
obtained~in~\cite{[Kro]}.%%%\enlargethispage{12pt}

\begin{prop}\label{prop:repr-Kron}\mbox{}
  \begin{enumerate}
    
  \item If $\alpha\notin\Delta_{+}$, then the set of indecomposable
    representations of $\Kron$ with the dimension $\alpha$ is
    empty.
    
  \item If $\alpha\!\in\!\Delta^{\mathrm{re}}_{+}$, then an
    indecomposable representation of $\Kron$ with the dimension
    $\alpha$ is either the representation $\rho_n$ if $\alpha=(n+1,n)$
    or $\bar{\rho}_n$ if $\alpha=(n,n+1)$, $n\in\oZ_+$,
    
  \item If $\alpha\in\Delta^{\mathrm{im}}_{+}$, then
    $\setI_{\alpha}=\oC\oP^1$ and indecomposable representations have
    the dimension $(n+1,n+1)$ for some $n\in\oZ_+$.
  \end{enumerate}
\end{prop}%%%\enlargethispage{\baselineskip}

\footnotesize
\addtolength{\parskip}{-8pt}

\end{document}